\documentclass[10pt]{amsart}

\usepackage{amsmath, amsthm}
\usepackage{eucal}

\usepackage{amssymb}
\usepackage{amscd}
\usepackage{palatino}
\usepackage{latexsym}
\usepackage{epsfig}
\usepackage{graphicx}
\usepackage{amsfonts}
\usepackage{psfrag}

\usepackage{stmaryrd}

\usepackage{url}
\usepackage{cite}
\usepackage{enumitem}

\usepackage[margin=1in]{geometry}

\input xy
\xyoption{all}
\UseComputerModernTips

\numberwithin{equation}{section}
\numberwithin{figure}{section}

\newtheorem{theorem}{Theorem}[section]
\newtheorem{lemma}[theorem]{Lemma}
\newtheorem{proposition}[theorem]{Proposition}

\newtheorem{remark}[theorem]{Remark}
\newtheorem{example}[theorem]{Example}

\theoremstyle{definition}
\newtheorem{definition}[theorem]{Definition}

\newcommand{\C}{{\mathbb{C}}}

\newcommand{\R}{{\mathbb{R}}}

\newcommand{\n}{\mathfrak{n}}
\renewcommand{\t}{\mathfrak{t}}
\newcommand{\g}{\mathfrak{g}}
\newcommand{\h}{\mathfrak{h}}

\renewcommand{\k}{\mathfrak{k}}
\renewcommand{\a}{\mathfrak{a}}

\newcommand{\into}{\hookrightarrow}

\renewcommand{\mod}{\mathbin{/\!/}}

\DeclareMathOperator{\im}{im}

\DeclareMathOperator{\Ad}{Ad}

\DeclareMathOperator{\tr}{Tr}

\begin{document}

\title{A local normal form for Hamiltonian actions of compact semisimple Poisson-Lie groups}

\author{Megumi Harada}
\address{Department of Mathematics and
Statistics\\ McMaster University\\ 1280 Main Street West\\ Hamilton, Ontario L8S4K1\\ Canada}
\email{Megumi.Harada@math.mcmaster.ca}
\urladdr{\url{http://www.math.mcmaster.ca/Megumi.Harada/}}

\author{Jeremy Lane}
\address{Department of Mathematics and
Statistics\\ McMaster University\\ 1280 Main Street West\\ Hamilton, Ontario L8S4K1\\ Canada}
\email{lanej5@math.mcmaster.ca}

\author{Aidan Patterson}
\address{Department of Mathematics\\ University of Waterloo \\ 200 University Avenue West \\ Waterloo,  Ontario\\ Canada}
\email{ae2patte@uwaterloo.ca}


\date{\today}


\begin{abstract}

The main contribution of this manuscript is a local normal form for Hamiltonian actions of Poisson-Lie groups $K$ on a symplectic manifold 
equipped with an $AN$-valued moment map, where $AN$ is the dual Poisson-Lie group of $K$. Our proof uses the delinearization theorem of Alekseev which 
relates a classical Hamiltonian action of $K$ with $\mathfrak{k}^*$-valued moment map to a Hamiltonian action with an $AN$-valued moment map, via a deformation  of symplectic structures.  We obtain our main result 
by 
proving a ``delinearization commutes with symplectic quotients'' theorem which is also of independent interest, and then putting this together with the local normal form theorem for classical Hamiltonian actions wtih $\k^*$-valued moment maps. 
A key ingredient for our main result is the delinearization $\mathcal{D}(\omega_{can})$ of the canonical symplectic structure on $T^*K$, so we additionally take some steps toward explicit computations of $\mathcal{D}(\omega_{can})$. In particular, in the case $K=SU(2)$, we 
obtain explicit formulas for the matrix coefficients of $\mathcal{D}(\omega_{can})$ with respect to a natural choice of coordinates on $T^*SU(2)$. 
\end{abstract}

\maketitle

\setcounter{tocdepth}{1}
\tableofcontents

\section{Introduction}\label{sec:intro}

Let $K$ be a compact, connected, semisimple Lie group and let $G$ be the complexification of $K$.  The Marle-Guillemin-Sternberg  local normal form for Hamiltonian $K$-actions on symplectic manifolds provides local models for such actions, and is an indispensable tool in Hamiltonian geometry \cite{GuilleminSternberg,Marle}. 
An Iwasawa decomposition $G = KAN$ endows $K$ with a Poisson-Lie group structure \cite{LuWeinstein}, and with respect to this structure, one can build an analogous theory of \emph{Hamiltonian Poisson-Lie group actions} of $K$ whose moment maps take values in the group $AN$ rather than the dual of the Lie algebra, $\mathfrak{k}^*$ \cite{LuMomentMap}. The main contribution of this article is a local normal form theorem, analogous to that of Marle-Guillemin-Sternberg, 
for Hamiltonian Poisson-Lie group actions with $AN$-valued moment maps. We expect that our local normal form can play a role similar to its classical counterpart in this setting. 

Our main tool is a theorem of Alekseev, which 
establishes a correspondence between classical Hamiltonian actions of $K$ with $\mathfrak{k}^*$-valued moment maps and Hamiltonian Poisson-Lie group actions of $K$ with $AN$-valued moment maps \cite{Alekseev}. This correspondence is achieved by a deformation of symplectic structures which works as follows. Let $M$ be a manifold equipped with a $K$-action. If $\omega$ is a symplectic form on $M$ with respect to which the $K$-action is a Hamiltonian Poisson-Lie group action with an $AN$-valued moment map, then there exists a deformed symplectic structure, $\mathcal{L}(\omega)$, on $M$ with respect to which the $K$-action is Hamiltonian in the classical sense. The deformed symplectic structure $\mathcal{L}(\omega)$ is called the \emph{linearization} of $\omega$. This result also works in reverse, and the reverse deformation is called \emph{delinearization}. The delinearization of a symplectic form is denoted $\mathcal{D}(\omega)$.

With this in mind, the basic idea of our manuscript is two-fold. We first prove a ``delinearization commutes with symplectic quotients" result (cf. Proposition~\ref{prop: delinearization quotient}) which is of independent interest, aligning as it does with similar well-known results in Hamiltonian geometry (such as the ``quantization commutes with symplectic quotients'' theorem).  We then use our ``delinearization commutes with quotients'' result, together with the Marle-Guillemin-Sternberg local normal form, to prove our local normal form theorem (cf. Theorem~\ref{theorem: local normal form for Lu mom maps}). More precisely, suppose that $(M,\omega)$ is a symplectic manifold equipped with a Hamiltonian Poisson-Lie group action of $K$ with an $AN$-valued moment map and suppose that $p\in M$ is sent to the group identity $e\in AN$ by the moment map $\Psi$. Our local normal form for a neighbourhood of an orbit of $p$ is 
\[
	(T^*K \times W, \mathcal{D}(\omega_{can}) \oplus \omega_W)\sslash_0 K_p
\]
where $(T^*K,\omega_{can})$ is the cotangent bundle of $K$ with its canonical symplectic structure, $K_p$ is the stabilizer subgroup of $p$, $W$ is the symplectic cross-section at $p$ for the action of $K$, and $\omega_W$ is the induced symplectic structure on $W$. The notation $\sslash_0 K_p$ denotes a symplectic quotient with respect to a classic Hamiltonian action of $K_p$. We restrict to the special case where $\Psi(p)=e$ for simplicity; we expect the general case will be a natural generalization for follow-up work. See Section~\ref{sec: local} for further details.

What makes this normal form complicated, in comparison to the Marle-Guillemin-Sternberg local normal form,  is the presence of the deformed symplectic structure $\mathcal{D}(\omega_{can})$.  To put our local normal form to use, it is evident that a concrete description of $\mathcal{D}(\omega_{can})$ is necessary. Therefore, we devote the last two sections of the paper to a more explicit description of $\mathcal{D}(\omega_{can})$.  In Section~\ref{sec: delin formula}, we use moment map equations to analyze the coefficients of the symplectic Poisson structure of $\mathcal{D}(\omega_{can})$ in the general case. We only solve for those coefficients which can be described using general properties of delinearization. Note that, in general, the delinearization $\mathcal{D}(\omega_{can})$ is not unique, so it is impossible to solve for all coefficients without specifying a specific delinearization.  However, in the special case in which $K$ has rank 1, the delinearization $\mathcal{D}(\omega_{can})$ is unique. Using this, we solve for the remaining coefficients of the symplectic Poisson structure of $\mathcal{D}(\omega_{can})$ for the special case where $K = SU(2)$ in Section \ref{subsec: pi for SU2}. Finally, in Section~\ref{sec: SU2}, we give a more direct derivation (without going through the associated Poisson structure) of $\mathcal{D}(\omega_{can})$ in the case $K = SU(2)$. This computation uses a slightly different (but natural) set of coordinates than that used in Section~\ref{sec: delin formula}. Specifically, in Section~\ref{sec: SU2} we use coordinates on $\mathfrak{su}(2)^{*}$ dual to the coordinates $\{t,x,y\}$ on $\mathfrak{su}(2)$, which are then identified with elements of $\textrm{Lie}(AN) = \mathfrak{a} \mathfrak{n}$, and we give explicit formulas for the matrix entries of $\mathcal{D}(\omega_{can})$ with respect to these coordinates. 

We close with some remarks about motivation, and possible directions for future work. First, and as mentioned above, it would be natural to extend our local normal form result to the general case when $\Psi(p)$ is not necessarily equal to the identity $e \in K$. Second, both the Poisson-Lie group structure on $K$ and the delinearization $\mathcal{D}(\omega)$ of the symplectic structure of a manifold $M$ equipped with a classical Hamiltonian $K$-action can be made to depend on a non-zero scalar $s$. It was recently shown in \cite{AlekseevHoffmanLaneLi}, that, generically, $\mathcal{D}(\omega)$ converges to a canonical constant symplectic structure, known as ``action-angle coordinates,'' in the limit $s \to \infty$.  We hope that our local normal form will aid in investigating the $s \to \infty$ limit non-generically.

\bigskip

\noindent \textbf{Acknowledgements.}   This manuscript is an outgrowth of a summer undergraduate research project conducted by the third author in two consecutive summers, 2019 and 2020, under the joint supervision of the first and second authors. This research was made possible by the support of the James Stewart Fund of the Department of Mathematics and Statistics at McMaster University and the NSERC USRA program. 
The first author is additionally supported by an NSERC Discovery Grant and a Canada Research Chair Award.

\section{Background and notation}\label{sec:background} 

The purpose of this section is to briefly recall some background and set notation. 

\subsection{Poisson and symplectic structures}\label{subsec: Poisson background} 

For more details on what follows we refer to \cite{CannasDaSilva}.  A {\bf Poisson structure} on a smooth manifold $M$ is a smooth bivector field $\pi \in \Gamma(\wedge^2TM)$ such that the associated pairing
\begin{equation}
	\{\cdot ,\cdot\} \colon C^\infty(M) \times C^\infty(M) \to \R, \quad \{f,g\} = \pi(df,dg)
\end{equation}
satisfies the Jacobi identity. The pairing $\{\cdot ,\cdot\}$ above is called the {\bf Poisson bracket} of $\pi$. Being a 
bivector field, $\pi$ may be paired against cotangent vectors; we notate this pairing as  
\begin{equation}\label{eq: pi sharp}
	\pi^\sharp\colon T^*M \to TM, \quad \alpha \mapsto \pi(\alpha,\cdot).
\end{equation}
A Poisson structure is {\bf symplectic} if $\pi^\sharp$ is a fiberwise isomorphism.  A {\bf symplectic structure} on $M$ is a closed 2-form $\omega\in \Omega^2(M)$ such that
\begin{equation}\label{eq: omega flat} 
	\omega^\flat\colon TM \to T^*M, \quad X \mapsto \omega(X,\cdot).
\end{equation}
is a fiberwise isomorphism. A symplectic structure $\omega$ determines a symplectic Poisson structure $\pi$ by the equation $\pi^\sharp =-( \omega^\flat)^{-1}$.

\subsection{Lie theory}\label{subsec: L background} 
We now recall some facts and notation related to Lie groups, and Poisson-Lie groups in particular.

Let $G$  be a connected complex semisimple Lie group and let  $H\subset G$ be a maximal complex torus. Denote their Lie algebras by $\g$ and $\h$ respectively. Let $R$ denote the set of roots of $G$. For each root $\alpha \in R$, let $\g_\alpha$ denote the corresponding root space.  Fix a choice of positive roots $R_+$. The subalgebra $\n= \bigoplus_{\alpha \in R_+}\g_\alpha$ is the Lie algebra of a maximal unipotent subgroup $N \subset G$.

Let $\kappa$ denote the Killing form of $\g$,  $\kappa(X,Y) := \tr(ad_X\circ ad_Y)$. 
For each $\alpha \in R_+$, fix $e_\alpha \in \g_\alpha$ and $e_{-\alpha} \in \g_{-\alpha}$ such that $\kappa(e_\alpha,e_{-\alpha}) = 1$ and let $h_\alpha := [e_\alpha,e_{-\alpha}]$. Let $\dagger$ denote the $\C$-antilinear involutive Lie algebra anti-homomorphism of $\g$ defined by $e_\alpha^\dagger = e_{-\alpha}$ and $h_\alpha^\dagger = h_\alpha$ for all $\alpha \in R$.  (Alternatively, one may describe $\dagger$ as the composition of the Cartan involution given by complex conjugation of $\g = \k^\C$ with the Lie algebra anti-homomorphism $X \mapsto -X$.) Let $\tau \colon G \to G$ denote the anti-holomorphic involutive anti-homomorphism such that $d\tau_e(X) = X^\dagger$. The compact real form of $G$ is the real Lie group $K$ formed by the set of $g$ such that $\tau(g) = g^{-1}$.
The Lie algebra of $K$ is the real Lie algebra $\k$ formed by the set of $X\in \g$ such that $X^\dagger = -X$. Equivalently, $\k$ is the real Lie subalgebra of $\g$ spanned by
\begin{equation}\label{standard elements}
    x_\alpha = \frac{1}{\sqrt{2}}(e_\alpha - e_{-\alpha}), \quad y_\alpha = \frac{\sqrt{-1}}{\sqrt{2}}(e_\alpha + e_{-\alpha}), \quad \sqrt{-1}h_\alpha, \quad \alpha \in R.
\end{equation}
Let $T$ denote the (real) maximal torus $T = H \cap K $ and let $\t $ denote its (real) Lie algebra. Let $A$ denote the real Lie subgroup of $G$ with Lie algebra $\a = \sqrt{-1}\t$. Let $\a\n$ denote the Lie algebra of the product (real) Lie subgroup $AN \subset G$. The pairing
\begin{equation}\label{eq: negative Killing}
(X, Y) := - \kappa(X, Y) \, \, \textup{ for } X, Y \in \k
\end{equation} 
is an invariant inner-product on $\k$. Let $\{t_i\}_{i\in [\ell]}$ denote a choice of orthonormal basis for $\t$ with respect to this inner product. Then
\begin{equation}\label{eq: orthonormal basis}
	\mathcal{B} := \{x_\alpha, y_\alpha, t_i\}_{i\in [\ell],\alpha\in R_+}
\end{equation}
 is an orthonormal basis of $\k$. The dual basis of the  dual vector space $\k^*$ is denoted $\mathcal{B}^* = \{x_\alpha^*, y_\alpha^*, t_i^*\}_{i\in [\ell],\alpha\in R_+}$.  (At various points in the paper, we will also use  $\mathcal{B}^*$ to denote a particular basis of $\a\n$. Although these two bases of $\k^*$ and $\a\n$ will be identified via a linear isomorphism, it will be helpful to keep the distinction between the two in mind.) 

\begin{example}\label{example: general linear group} 
Let $G= SL_n(\C)$ and let $H$ be the maximal torus consisting of diagonal matrices in $G$. Then $\g = \mathfrak{sl}_n(\C)$ is the Lie algebra of traceless $n\times n$ complex matrices and $\h$ is the set of traceless diagonal matrices. The Killing form on $\g$ is $\kappa(X,Y) = 2n\tr(XY)$. The subgroup $A$ consists of diagonal matrices with positive real diagonal entries. Let $E_{i,j}$ denote the elementary $n\times n$ matrix with $(i,j)$-th entry equal to 1 and all other entries equal to 0. We fix $R_+$ so that the maximal unipotent subgroup $N$ is the subgroup of upper triangular matrices with diagonal entries equal to 1. With this choice, the set of positive roots is indexed by pairs $(i,j)$, $1 \leq i < j \leq n$ and the root subspace corresponding to $(i,j)$ is spanned by $E_{i,j}$. The anti-involution $\dagger$ is conjugate-transpose, i.e. $X^\dagger = \overline{X}^T$, and $\tau(g) = \overline{g}^T$. For $\alpha \in R_+$, the Lie algebra elements $x_\alpha, y_\alpha, h_\alpha$ are given by 
\begin{equation}\label{sln standard elements}
    x_{i,j} = \frac{1}{2\sqrt{n}}(E_{i,j} - E_{j,i}), \quad y_{i,j} = \frac{\sqrt{-1}}{2\sqrt{n}}(E_{i,j} + E_{j,i}), \quad h_{i,j} = \frac{E_{i,i} - E_{j,j}}{2n}, \quad 1 \leq i < j \leq n.
\end{equation}
In the special case $n=2$, there is a single positive root and the elements of our orthonormal basis $\mathcal{B}$ are:
\begin{equation}\label{sl2 orthonormal basis}
    x = \frac{1}{2\sqrt{2}}\left(\begin{array}{cc}
		0 & 1 \\ -1 & 0 
	\end{array}\right), \quad y = \frac{\sqrt{-1}}{2\sqrt{2}}\left(\begin{array}{cc}
		0 & 1 \\ 1 & 0 
	\end{array}\right), \quad t = \sqrt{-2}h = \frac{\sqrt{-1}}{2\sqrt{2}}\left(\begin{array}{cc}
		1 & 0 \\ 0 & -1 
	\end{array}\right).
\end{equation}
\end{example}

\begin{remark}\label{remark: semisimple}  
Throughout the remainder of this manuscript, we will assume for simplicity that the Lie group $K$ under consideration is semisimple (and hence fits into the framework discussed above). 
\end{remark} 

\subsection{Poisson-Lie groups}\label{subsec: PL background} 
We continue with the notation of Section~\ref{subsec: L background} and now briefly recall the theory of Poisson-Lie groups and their duals. See \cite{LuWeinstein} for more details.

A \textbf{Poisson-Lie group} is a Lie group equipped with a Poisson structure such that multiplication is a Poisson map. For a real Lie group $K$ arising in the setup of Section~\ref{subsec: L background}, there is a family of Poisson-Lie group structures on $K$ defined as follows. For all $s\neq 0$, the pairing 
\begin{equation}
	\langle\cdot,\cdot\rangle_s := \frac{-1}{s}\im \kappa(\cdot,\cdot) 
\end{equation}
is an invariant, non-degenerate, symmetric bilinear form on $\g$ such that the Lie subalgebras $\k$ and $\a\n$ are isotropic. In other words, the tuple $(\g,\k,\a\n,\langle\cdot,\cdot\rangle_s)$ is a Manin triple. This Manin triple determines Lie bialgebras $(\k,[,],\delta_{\k,s})$ and $(\a\n, [,], \delta_{\a\n,s})$ that are dual to one another with respect to the pairing $\langle\cdot,\cdot\rangle_s$. Note that whereas the Lie algebra structures on $\k$ and $\a\n$ are fixed, the Lie algebra 1-cocycles $\delta_{\k,s}$ and $\delta_{\a\n,s}$ are $s$-dependent because our pairing is $s$-dependent.

The  Lie bialgebra $(\k,[,],\delta_{\k,s})$ is a coboundary Lie bialgebra, i.e., there is a skew-symmetric $r$-matrix $r_{K,s} \in \wedge^2 \k$ such that $\delta_{\k,s} = \partial r_{K,s}$.  This $r$-matrix defines a Poisson-Lie group structure on $K$ by the formula
\begin{equation}\label{eq: pi K s}
	\pi_{K,s} = r_{K,s}^L - r_{K,s}^R
\end{equation}
where $r_{K,s}^L$ (respectively, $r_{K,s}^R$) denotes the bivector field on $K$ that is invariant with respect to left (respectively, right) multiplication and satisfies $r_{K,s}^L(e) = r_{K,s}$ (respectively, $r_{K,s}^R(e) = r_{K,s}$).  More explicitly, if $L_k$ denotes left multiplication by $k$ and $R_k$ denotes right multiplication by $k^{-1}$, then $r_{K,s}^L(k) = (L_k)_*r_{K,s}$ and $r_{K,s}^R(k) = (R_{k^{-1}})_*r_{K,s}$.  The tangent Lie bialgebra of $(K,\pi_{K,s})$ is $(\k,[,],\delta_{\k,s})$. The parameter $s$ corresponds to a linear scaling, with $\pi_{K,s} = s\pi_{K,1}$.  Thus setting $s=0$ in this last equation recovers the trivial Poisson-Lie group $(K,0)$. 
\begin{example} 
In the case $\k = \mathfrak{su}(2)$, the $r$-matrix is $r_{K,s} = s \cdot x\wedge y$.
\end{example}

There is also a Poisson-Lie group structure $\pi_{AN,s}$ on $AN$ whose tangent Lie bialgebra is $(\a\n, [,], \delta_{\a\n,s})$. The Poisson-Lie groups $(K,\pi_{K,s})$ and $(AN,\pi_{AN,s})$ are said to be \textbf{dual} to one another because their tangent Lie bialgebras are dual to one another. We now recall an explicit expression for $\pi_{AN,s}$ that will be useful in what follows. We refer the reader to 
\cite[Chapter 11]{PoissonStructures}, in particular \cite[Proof of Theorem 11.47]{PoissonStructures}, for more details (note that the definition of ``dual Poisson-Lie group'' given there differs from the definition used here by a sign). 
 
The \textbf{double} of the dual Lie bialgebras $(\k,[,],\delta_{\k,s})$ and $(\a\n, [,], \delta_{\a\n,s})$ is the (real) Lie algebra $\g = \k \oplus \a\n$. Let $\mathcal{B}$ be the orthonormal basis of $\k$ which was fixed above.  Along with the non-degenerate pairing  $\langle\cdot,\cdot\rangle_s$, this determines a dual basis $\mathcal{B}^* \subset \a\n$.  Note that since the pairing is $s$-dependent, the dual basis $\mathcal{B}^*$ is also $s$-dependent. A skew-symmetric $r$-matrix for the double can be expressed in terms of these dual bases as
\[
r_{G,s} = \frac{1}{2} \sum_{b\in \mathcal{B}} b \wedge b^*.
\]
This $r$-matrix defines a Poisson-Lie group structure $\pi_{G,s} = r_{G,s}^L - r_{G,s}^R$ on the double $G$ with the property that $\pi_{AN,s} = - \pi_{G,s} \vert_{AN}$ and $\pi_{K,s} = \pi_{G,s} \vert_{K}$. It follows that 
\begin{equation}\label{eq: piK star}
\pi_{AN,s} =( r_{G,s}^R - r_{G,s}^L) \bigg\vert_{AN} \quad \text{and} \quad \pi_{K,s} =( r_{G,s}^L - r_{G,s}^R) \bigg\vert_{K}.
\end{equation} 
Let $Pr_{\mathfrak{an}}: \mathfrak{g} \to \mathfrak{an}$  denote projection with respect to the decomposition $\mathfrak{g} = \mathfrak{k} \oplus \mathfrak{an}$. 

\begin{lemma}\label{lemma: piAN take 1} 
	Let $\mathcal{B}$ be the  orthonormal basis \eqref{eq: orthonormal basis} of $\mathfrak{k}$, and let $\mathcal{B}^*$ be the 
	dual basis of $\mathfrak{an}$ determined by the pairing  $\langle\cdot,\cdot\rangle_s$. Then, for all $p \in AN$,
	\begin{align}\label{PoissonBivectorAN}
	(L_{p^{-1}})_*(\pi_{AN,s})_{p} &= \frac{1}{2}\sum_{b\in \mathcal{B}}  Pr_{\mathfrak{an}}(Ad_{p^{-1}}(b))\wedge Ad_{p^{-1}}(b^*). 
	\end{align} 
\end{lemma}

\begin{proof}
	We follow the argument in the proof of \cite[Theorem 11.47]{PoissonStructures}. 
	By~\eqref{eq: piK star},
	\[
	\pi_{AN,s} =  ( r_{G,s}^R - r_{G,s}^L) \bigg\vert_{AN} = \left(\frac{1}{2} \sum_{b\in \mathcal{B}} b \wedge b^* \right)^R\bigg\vert_{AN}  - \left( \frac{1}{2} \sum_{b\in \mathcal{B}} b \wedge b^* \right)^L \bigg\vert_{AN} .
	\]
	It follows that for $p \in AN$, 
	\begin{align*}
	(L_{p^{-1}})_*(\pi_{AN,s})_p &= \frac{1}{2} \left(\sum_{b\in \mathcal{B}} Ad_{p^{-1}}(b)\wedge Ad_{p^{-1}}(b^*) - \frac{1}{2}\sum_{b\in \mathcal{B}} b\wedge b^*\right)  \bigg\vert_{\a\n}. \\
	\end{align*}
	The vectors $b \in \mathcal{B}$ pair trivially against element in $T^*_e(AN)$. We also observe that, since $p \in AN$ and $b^* \in \mathfrak{an}$, the vector $Ad_{p^{-1}}(b^*)$ is automatically an element of $\mathfrak{an}$. The formula follows.
\end{proof}

The Iwasawa decompositions of $G$ define diffeomorphisms: $K\times AN \cong G$, $(k,b) \mapsto kb$, and $AN \times K \cong G$, $(b,k) \mapsto bk$. Let $pr_{K}\colon G \to K$ and $pr_{AN} \colon G \to AN$ denote the projections to the right and left factors, respectively, with respect to the decomposition $G\cong AN \times K$.  The \textbf{dressing action} of $K$ on $AN$ is defined to be the left action
\begin{equation}\label{left dressing action}
	K \times AN \to AN, \quad (k,b) \mapsto b^k := pr_{AN}(kb).
\end{equation}
The dressing action of $AN$ on $K$ is similarly defined as the right action given by 
\begin{equation}\label{right dressing action}
	K \times AN \to K, \quad (k,b) \mapsto k^b := pr_{K}(kb).
\end{equation}

We now describe an important map $E_s \colon \k^* \to AN$, which depends on $s$, which will be a key ingredient in our computations below. More precisely, for all $s\neq 0$, define
\begin{equation}\label{Es definition}
	E_s \colon \k^* \xrightarrow{\phi }\k \xrightarrow{\exp(2s\sqrt{-1} \cdot ) } P \xrightarrow{f^{-1}} AN
\end{equation}
where: $\phi$ is the linear isomorphism defined by the inner-product~\eqref{eq: negative Killing}, the space $P$ is defined to be the image $ \exp(\sqrt{-1}\k)\subset G$ of $\sqrt{-1}\k$ under the exponential map, and $f^{-1}$ is the inverse of the diffeomorphism
\begin{equation}\label{eq: def f Kstar to P} 
	f\colon AN \to P,\quad b \mapsto b\tau(b).
\end{equation}
It is useful to notice that $E_s$ has a simple form when restricted to the dual of the torus of $K$. This is because 
for $\lambda \in \t^*$ it can be computed that the composition $\t^* \xrightarrow{\phi} \k \xrightarrow{\exp(2s\sqrt{-1} \cdot)} P$ has image landing in $A$, the subgroup of diagonal matrices with real positive entries. The inverse $f^{-1}$ of $f$ defined in~\eqref{eq: def f Kstar to P}, when restricted to $A$, is the map which simply takes the square root of each entry. From this it is possible to see that for $\lambda \in \t^*$ we have 
\begin{equation}\label{eq: Es on torus} 
E_s(\lambda) = \exp(s \sqrt{-1} \phi(\lambda)).
\end{equation} 
The map $E_s: \k^* \to AN$ is a diffeomorphism and it is equivariant with respect to the coadjoint action of $K$ on $\k^*$ and the dressing action of $K$ on $AN$ \cite{AlekseevMeinrekenWoodward}.  Thus, for all $k\in K$ and $\lambda \in \k^*$ such that $\phi(\lambda) \in\t$, we have $E_s(Ad_k^*\lambda) = E_s(\lambda)^k = \exp(s\sqrt{-1}\phi(\lambda))^k$ where the $k$ in the superscript denotes the left dressing action of~\eqref{left dressing action} and we have used~\eqref{eq: Es on torus} above. 

\subsection{Hamiltonian actions of Poisson-Lie groups}\label{subsec: HPL background} 

Let $K$ be a compact connected semisimple (real) Lie group. 
Let $(M,\omega)$ be a symplectic manifold.  Given a smooth left action of $K$ on $M$, the \textbf{infinitesimal generated vector field} of $X\in \k$ on $M$ is the smooth vector field $X_M$ given by 
\begin{equation}\label{eq: XM} 
	X_M(m) := \frac{d}{dt}\bigg\vert_{t=0} \exp(tX)\cdot m, \quad m \in M.
\end{equation}
Let $\Theta^L$ denote the left-invariant Maurer-Cartan form on $AN$. 

\begin{definition}\label{definition: moment map} 
A {\bf Hamiltonian action of $(K,\pi_{K,s})$ on $(M,\omega)$}  is an action of $K$ on $M$ such that there exists a $K$-equivariant map $\Psi: M \to AN$ with  
\begin{equation}\label{K* valued moment map}
	\omega(X_M,\cdot) = \Psi^*\langle \Theta^L,X\rangle_s \quad \forall X\in \k.
\end{equation}
A {\bf Hamiltonian action of $(K,0)$ on $(M,\omega)$}  is an action of $K$ on $M$ such that there exists a $K$-equivariant map $\mu: M \to \k^*$ with 
\begin{equation}\label{k* valued moment map}
	\omega(X_M,\cdot) = d\langle \mu,X\rangle  \quad \forall X\in \k.
\end{equation}
\end{definition}  
Maps satisfying \eqref{K* valued moment map} or \eqref{k* valued moment map} are called {\bf moment maps}. 
Although the moment map equations \eqref{K* valued moment map} and \eqref{k* valued moment map} are quite similar, there are also significant differences. For instance, Hamiltonian actions of $(K,0)$ preserve the symplectic structure, while Hamiltonian actions of $(K,\pi_{K,s})$ do not. Instead, the action map $\sigma\colon K \times M \to M$, $\sigma(k,p) = k\cdot p$ is a Poisson map with respect to the product Poisson structure on $K\times M$ and the symplectic Poisson structure of $M$. Equivalently, we have \cite[Theorem 2.6]{LuWeinstein}
\begin{equation}\label{action map is Poisson}
	\pi(k\cdot p) = \sigma(k)_* \pi(p) + \sigma(p)_* \pi_{K,s}(k).
\end{equation}

Despite these differences, Hamiltonian actions of $(K,0)$ and $(K,\pi_{K,s})$ are equivalent in a specific sense, described by the theorem below.  We need some preliminaries. Suppose $\Omega^s \in \Omega^2(\k^*)$ is a closed  $2$-form such that 
	\begin{equation}\label{eq: Omega s def}
		\Omega^s(X_{\k^*}, \cdot ) = E_s^* \langle \Theta^L, X \rangle_s  -  d \langle \cdot, X \rangle \quad \forall X\in \k
	\end{equation}
	where $X_{\k^*}$ is the infinitesimal generated vector field of $X\in \k$ on $\k^*$ with respect to the coadjoint action of $K$ on $\k^*$. There exist  $2$-forms $\Omega^s$ with this property, but  this property does not determine $\Omega^s$ uniquely  \cite{Alekseev, AlekseevMeinrekenWoodward}. We have the following. 
	
\begin{theorem}\label{theorem: delinearization} \cite{Alekseev}
Let $K$ be as above and let $K \times M \to M$ be a smooth left action of $K$  on $M$. Let $\Omega^s \in \Omega^2(\k^*)$ be a  closed $2$-form such that \eqref{eq: Omega s def} holds.  Suppose $\Psi\colon M \to AN $ and $\mu\colon M \to \k^*$ are smooth $K$-equivariant maps related by 
	\begin{equation}\label{eq: delin 1} 
		\Psi = E_s \circ \mu. 
	\end{equation} 
	Suppose in addition that $\omega,\omega'\in \Omega^2(M)$ are two symplectic structures related by 
	\begin{equation}\label{eq: delin 2} 
		\omega' = \omega +  \mu^* \Omega^s.
	\end{equation} 
	Then the action of  $(K,\pi_{K,s})$ on $(M,\omega')$  is Hamiltonian with moment map $\Psi$ if and only if the action of  $(K,0)$ on $(M,\omega)$  is Hamiltonian with moment map $\mu$.
\end{theorem}

If $\omega'$ and $\omega$ are  related as  in~\eqref{eq: delin 2} then we say that $\omega'$ is the \textbf{delinearization} of $\omega$ and $\omega$ is the  \textbf{linearization} of $\omega'$. 
Given a fixed choice of $\Omega^s$ satisfying \eqref{eq: delin 2},  we denote delinearization and linearization as 
\begin{equation}
	\mathcal{D}(\omega) := \omega +  \mu^* \Omega^s, \quad\text{and}\quad \mathcal{L}(\omega) := \omega - \mu^* \Omega^s
\end{equation}
respectively.  Note that $\mathcal{D}$ and $\mathcal{L}$ both depend on the parameter $s$, although we suppress this dependence from the notation. The 2-forms $\Omega^s$ can be constructed to be exact, so $[\omega] = [\mathcal{D}(\omega)]$ in cohomology. The next lemma will be useful in what follows.
\begin{lemma}\label{background useful fact}
Let $(M,\omega)$ be a symplectic manifold equipped with a Hamiltonian action of $(K,\pi_{K,s})$ with moment map $\Psi\colon M \to AN$. Then, for all $X,Y \in \k$ and $m \in M$,
\begin{align}
	\omega_m(X_M,Y_M)   = (\pi_{AN,s})_p( \iota^*\langle\Theta^L,X\rangle_s,\iota^*\langle\Theta^L,Y\rangle_s) =  (L_{p^{-1}})_*(\pi_{AN,s})_p(X,Y)
\end{align}
where $p = \Psi(m)$ and in the last expression $X$ denotes $ \langle -, X\rangle_s \in \mathfrak{an}^*$,  and $Y$ denotes $ \langle -, Y\rangle_s \in \mathfrak{an}^*$.
\end{lemma}

\begin{proof}
We have 
\begin{align*}
	\omega_m(X_M,Y_M) & = (\Psi^*\langle\Theta^L,X\rangle_s)_m(Y_M) \\
	& = (\langle\Theta^L,X\rangle_s)_p(\Psi_*Y_M)\\
	& = ( \langle\Theta^L,X\rangle_s)_p(Y_{AN})\\
	& = (\iota^*\langle\Theta^L,X\rangle_s)_p (Y_{AN})\\
	& = (\pi_{AN,s})_p(\iota^*\langle\Theta^L,X\rangle_s,\iota^*\langle\Theta^L,Y\rangle_s)
\end{align*}
where the first equality is the moment map equation \eqref{K* valued moment map}. The third equality is $K$-equivariance of the moment map $\Psi$, where $Y_{AN}$ denotes the dressing vector field of $Y$ on $K^*$.  Since the dressing vector field is tangent to the dressing orbits, this formula is unchanged if we replace the 1-form $\langle\Theta^L,X\rangle_s$ with its pullback $\iota^*\langle\Theta^L,X\rangle_s$ to the dressing orbit through $\Psi(p)$ for whichever point $p$ we are evaluating this real-valued function at (hence the fourth equality). The fifth equality follows since the dressing action of $K$ on $K^*$ restricts to a Hamiltonian action of $(K,\pi_{K,s})$ on the dressing orbit through $\Psi(p)$, with moment map given by the inclusion map $\iota$. Finally, note that $\langle \Theta^L, X\rangle_s(p) = L_{p^{-1}}^*X$ so
\[
		(\pi_{AN,s})_p(\iota^*\langle \Theta^L,X\rangle_s,\iota^*\langle \Theta^L,Y\rangle_s)  = (\pi_{AN,s})_p(L_{p^{-1}}^*X,L_{p^{-1}}^*Y)  = (L_{p^{-1}})_*(\pi_{AN,s})_p(X,Y). \qedhere
\]
\end{proof}

\subsection{Marle-Guillemin-Sternberg local normal forms}\label{subsec: local normal forms background} 

In this section we recall the Marle-Guillemin-Sternberg local normal form theorem. Let $K$ (or $H$) denote a compact connected semisimple (real) Lie group.  Throughout this section, all Lie groups have trivial Poisson structure, so we drop the Poisson-Lie notation. In particular, a \textbf{Hamiltonian $K$-action} is a Hamiltonian $(K,0)$-action in the sense of Definition \ref{definition: moment map}.

We first recall notation concerning Marsden-Weinstein symplectic quotients (cf. \cite{CannasDaSilva}). Suppose there is a Hamiltonian $H$-action on $(M, \omega)$ 
with moment map $\nu$.  Assuming $H$ acts freely on the level set $\nu^{-1}(0)$, we denote by 
$(M,\omega) \mod_0 H$ the \textbf{symplectic quotient} with underlying manifold $\mu^{-1}(0)/H$ and symplectic structure $\widetilde{\omega}$ defined by $pr^*\tilde{\omega} = \iota^*\omega$
where $\iota \colon \nu^{-1}(0) \to M$ is inclusion and $pr \colon \nu^{-1}(0) \to \nu^{-1}(0)/H$ is the quotient map. 
We sometimes denote the symplectic quotient as simply $M \mod_0 H$.

Next we recall the canonical symplectic structure on the cotangent bundle of a Lie group. Let $K$ be a compact connected Lie group,  $T^*K$ its cotangent bundle, and $\theta_{taut}$ the canonical tautological $1$-form on any cotangent bundle (cf. \cite{CannasDaSilva}). We take the convention that the \textbf{canonical symplectic form $\omega_{can}$} on $T^*K$ is $\omega_{can} = - d \theta_{taut}$. 
Fix the trivialization of the cotangent bundle $T^*K \cong K \times \k^*$ by  left-invariant 1-forms via the correspondence $\xi^L_k \leftrightarrow (k,\xi)$.  With this identification, the cotangent bundle $T^*(T^*K) \cong T^*K \times T^*\k^*$ has the frame 
\begin{equation}\label{eq; cotangent frame}
	\{ (a^*)^L \mid a^* \in \mathcal{B}^*\} \cup \{ b \mid  b \in \mathcal{B}\}
\end{equation}
where $(a^*)^L$ denotes the left invariant 1-form on $K$ with $(a^*)^L_e = a^*$, and $b$ is the constant vector field defined by the canonical identification $T^*\k^* \cong \k^* \times \k$.  This determines a frame
\begin{equation}\label{eq; wedge 2 cotangent frame}
	\{ (a^*)^L\wedge (b^*)^L \mid a^*,b^* \in \mathcal{B}^*\} \cup \{ a\wedge (b^*)^L \mid  a \in \mathcal{B},\, b^* \in \mathcal{B}^*\} \cup \{ a\wedge b \mid  a,b \in \mathcal{B}\}
\end{equation}
of $\wedge^2 T^*(T^*K)$. With respect to this frame,
\begin{equation}\label{omega can explicit}
\omega_{can} =  -\sum_{b\in \mathcal{B}} b \wedge (b^*)^L  +\frac{1}{2} \sum_{a,b \in \mathcal{B}} [a,b] (a^*)^L\wedge (b^*)^L
\end{equation}
where $[a,b]$ denotes the function $(k,\xi) \mapsto \langle \xi,[a,b]\rangle$.
\begin{example} In the case $K=SU(2)$ we have 
\[
	\omega_{can} = -t \wedge (t^*)^L - x\wedge (x^*)^L - y \wedge (y^*)^L + \frac{1}{\sqrt{2}}\left(t (x^*)^L \wedge(y^*)^L + y (t^*)^L \wedge (x^*)^L + x(y^*)^L \wedge (t^*)^L\right). 
\]
\end{example}

Next we recall two standard Hamiltonian structures on $T^*K$. With respect to the trivialization $T^*K \cong K \times \k^*$, the cotangent lift of the standard $K\times K$-action (by left and right multiplication) on $K$ is 
\begin{equation}\label{eq: K by K action on T*K} 
	L_{k_1}R_{k_2}(k,\xi) = (k_1kk_2^{-1},Ad_{k_2}^*\xi).
\end{equation}
The infinitesimal generated vector field of $(X,Y) \in \k \times \k$ is 
\begin{equation}
	(X,Y)_{T^*K}(k,\xi) = (X^R_k-Y^L_k,ad_Y\xi)  \in T_kK \times T_\xi\k^*.
\end{equation}
It is useful to note that $X_k^L = (Ad_kX)^R_k$. 
This $K\times K$ action on $(T^*K,\omega_{can})$ is Hamiltonian with moment map
\begin{equation}\label{eq: muL and muR} 
	(\mu_L,\mu_R) \colon K \times \k^* \to \k^* \times \k^*, \quad \mu_L(k,\xi) = Ad_k^*\xi,  \quad \mu_R(k,\xi) = -\xi.
\end{equation}

The following computation will be used below; the proof is straightforward. 

\begin{lemma}\label{lemma: dmu} 
Let $(X^R_k, \eta)$ be a tangent vector in $T_{(k,\lambda)}(T^*K) \cong T_kK \times \k^*$. Then 
	\begin{equation}\label{d mu}
	(d\mu_L)_{(k,\lambda)} \left( X^R_k,\eta \right) = ad^*_X(Ad_k^*(\lambda)) + Ad^*_k(\eta). 
	\end{equation}
\end{lemma}

Finally, we define the data we will need to express the local normal form theorems in what follows. Let $(M,\omega)$ be a symplectic manifold equipped with a Hamiltonian $K$-action with moment map $\mu: M \to \k^*$. To a point $p \in M$, we associate three pieces of {\bf local normal form data}: the stabilizer subgroup $K_p \subset K$ of the point $p$ with Lie algebra denoted $\k_p$; the moment map image $\mu(p) \in \k^*$; and the {\bf symplectic cross-section} at $p$, defined as 
\[
        W := T_p(K\cdot p)^\omega /( T_p(K\cdot p)^\omega \cap  T_p(K\cdot p)),
\]
where $T_p(K\cdot p)^\omega \subset T_pM$ denotes the symplectic complement of the subspace $T_p(K\cdot p)$ with respect to $\omega$. The symplectic cross-section has an induced linear symplectic structure denoted $\omega_W$, there is an induced linear action of $K_p$ on $W_p$ called the {\bf isotropy representation} which is Hamiltonian with respect to $\omega_W$, with quadratic moment map $\Phi_W \colon W \to \k_p^*$ defined by
\begin{equation}
	\frac{1}{2}\omega((X_W)_{v},v) = \langle\Phi_W(v),X\rangle, \quad \forall v\in W,\, X \in \k_p. 
\end{equation}
Note that $\Phi_W(0) = 0$. 
    
We can now recall the local model for the local normal form theorem.  Let $(K_p, \xi = \mu(p), W)$ be the local normal form data at $p$ as above. For simplicity we assume that $\xi = 0$. Let $K_p$ act diagonally on $T^*K \times W$ by
\begin{equation}
h \cdot (k,\eta, w) = (kh^{-1},\Ad_{h}^*\eta,h\cdot w)
\end{equation}
It follows from the above discussion that this $K_p$-action is Hamiltonian (with respect to the product symplectic structure) with moment map
\begin{equation}
	K\times \k^* \times W \to \k^*_p, \quad  (k,\eta,w) \mapsto -\eta\vert_{\k_p} +\Phi_W(w)
\end{equation}
where $-\eta\vert_{\k_p}$ denotes the image of $-\eta$ under the natural projection map $\k^* \to \k_p^*$. 
The {\bf local model} for data $K_p$, $0= \mu(p)$, $W$ is the symplectic quotient 
\begin{equation}\label{eq: local model space}
	(T^*K \times W,\omega_{can}\oplus \omega_W)\mod_0 K_p.
\end{equation}
Since the action of $K_p$ on the $0$-level set in $T^*K \times W$ is free, the space~\eqref{eq: local model space} is a smooth symplectic manifold. It is equipped with a Hamiltonian $K$-action which descends from the left action of $K$ on $T^*K$ defined by $L_{k'}(k,\eta) = (k'k,\eta)$. This action of $K$ has a moment map $\widetilde\mu_L$ which descends from the moment map $\mu_L$,
\begin{equation}
	\widetilde\mu_L \colon (T^*K \times W,\omega_{can}\oplus \omega_W)\mod_0 K_p \to \k^*, \quad \widetilde\mu_L([k,\eta, w]) = Ad_k^*(\eta)
\end{equation}
where $[k,\eta, w]$ denotes an equivalence class in the 0-level set under the action of $K_p$.

We have the following \cite{Marle,GuilleminSternberg}. 

\begin{theorem}[Marle-Guillemin-Sternberg]\label{theorem: local normal form}
Let $(M,\omega)$ be a symplectic manifold equipped with a Hamiltonian $K$-action with moment map $\mu: M \to \k^*$. Let $p \in \mu^{-1}(0)$. Then 
there exists a $K$-equivariant symplectomorphism $\varphi$ from a $K$-invariant neighbourhood of  $[e,0,0] $ in $(T^*K \times W,\omega_{can}\oplus \omega_W)\mod_0 K_p$ to a $K$-invariant neighbourhood of $p$ in $(M,\omega)$ such that $\varphi([e,0,0]) = p$ and $\mu\circ \varphi = \widetilde\mu_L$.
\end{theorem}

\section{Delinearization commutes with symplectic quotient}\label{sec: commutes} 

The main result of this section shows that the delinearization procedure outlined in Section~\ref{subsec: HPL background} commutes, in certain situations, with the Marsden-Weinstein symplectic quotient procedure described in Section~\ref{subsec: local normal forms background}.  
 This will be a key step in the proof of our main theorem in Section~\ref{sec: local}, but it is also of interest in its own right.

Let $H$ and $K$ be compact connected semisimple Lie groups and let $(M, \omega)$ be a symplectic manifold. Assume we are given a Hamiltonian action of $(H \times K,0)$ on  $(M, \omega)$ with moment map $(\nu, \mu): M \to \h^* \times \k^*$. Note in particular that this means that the actions of $H$ and $K$ commute and $\nu$ and $\mu$ can be separately viewed as moment maps for the actions of $H$ and $K$ respectively. Moreover, since the coadjoint action of $H$ on the second component of $\h^* \times \k^*$ is trivial (and similarly vice versa), we see that $\mu$ is $H$-invariant and $\nu$ is $K$-invariant. 
Finally, we assume for simplicity that the $H$-action on $M$ 
is \emph{free}. 
It is worth emphasizing here that the roles played by the $H$ and $K$ in the discussion below are quite different. Specifically, we will be ``delinearizing with respect to $K$'' and ``quotienting with respect to $H$''.

We begin with the following simple lemma.

\begin{lemma}\label{lemma: delinearized Hamiltonian by H}
Let 
$\mathcal{D}(\omega)$ denote the delinearization of the symplectic form $\omega$ with respect to the action of $K$, its
moment map $\mu$, a choice of parameter $s$, and corresponding $2$-form $\Omega^s$. Then the action of $H$ on $(M, \mathcal{D}(\omega))$ is Hamiltonian with moment map $\nu: M \to \h^*$. 
\end{lemma}

\begin{proof} 
Since we already know that 
$\nu: M \to \h^*$ is $H$-equivariant, it suffices 
 to show that $\nu$ satisfies the moment map equation~\eqref{k* valued moment map} with respect to the action of $H$ on $(M,\mathcal{D}(\omega))$. 
Since $\nu$ is a moment map with respect to the original symplectic structure $\omega$, we have from~\eqref{k* valued moment map} that 
	\begin{equation}\label{eq: mom map eqn for omega}
	\omega(Y_M, -) = d \langle \nu, Y \rangle 
	\end{equation} 
	where $Y \in \h$ and $Y_M$ denotes the infinitesimal generated vector field of $Y$ on $M$ with respect to the $H$-action. 
	Therefore, since $\mathcal{D}(\omega) = \omega + \mu^* \Omega^s$ it suffices to show that 
	\[
	(\mu^* \Omega^s)(Y_M, -) = 0
	\]
	for all $Y \in \h$. We first compute 
	\[
	(\mu^* \Omega^s)(Y_M, -)  = \mu^* (\Omega^s(d\mu(Y_M), \cdot)). 
	\]
	Now it would suffice to show that $d\mu(Y_M) = 0$ for all $Y \in \h$, but this follows from the fact that the $K$-moment map $\mu$ is invariant with respect to the action of $H$, by assumption. 
\end{proof}

From Lemma~\ref{lemma: delinearized Hamiltonian by H} it follows that the symplectic quotient $(M,\omega) \mod_0 H$ can be equipped with two natural delinearized symplectic structures. On the one hand, we can ``first delinearize, then take a symplectic quotient''.  Specifically, we may first consider the delinearization $\mathcal{D}(\omega)$ with respect to the $K$-action and the moment map $\mu$. By Lemma~\ref{lemma: delinearized Hamiltonian by H}, the group $H$ still acts in a Hamiltonian fashion on $(M, \mathcal{D}(\omega))$ with moment map $\nu$. Hence we may define a reduced symplectic structure $\widetilde{\mathcal{D}(\omega)}$ on $(M, \mathcal{D}(\omega)) \mod_0 H$ by the property that
\begin{equation}
	pr^* \widetilde{\mathcal{D}(\omega)} = \iota^*\mathcal{D}(\omega)
\end{equation}
where $pr\colon \nu^{-1}(0) \to \nu^{-1}(0)/H$ is the quotient projection map and $\iota \colon \nu^{-1}(0) \to M$ is the inclusion map.

On the other hand, we can also ``first take symplectic quotient, then delinearize''.  More precisely, we may first take 
$\widetilde{\omega}$ to be the quotient symplectic structure on $(M,\omega)\mod_0 H$.
It is well-known (see e.g. \cite{CannasDaSilva}) that, since the actions of $K$ and $H$ commute, 
the action of $K$ and the moment map $\mu$ descends to a 
Hamiltonian action on $(M,\omega)\mod_0 H$ with a moment map $\widetilde{\mu}$ naturally induced by $\mu$. 
Delinearizing $\widetilde{\omega}$ with respect to this induced Hamitonian $K$-action and moment map $\widetilde{\mu}$ results in a symplectic form 
\begin{equation}
\mathcal{D}(\widetilde{\omega}) = \widetilde{\omega} + \widetilde \mu^* \Omega^s.
\end{equation}

The main result of this section, below, states that these two constructions yield the same result. 
\begin{proposition}\label{prop: delinearization quotient} 
	(``Delinearization commutes with symplectic quotient''.)  Let $(M,\omega)$ be a symplectic manifold equipped with a Hamiltonian $(H \times K, 0)$-action with moment map $(\nu, \mu): M \to \h^* \times \k^*$. Assume that the $H$-action on $M$ is free. Let $s>0$ be a positive real parameter and suppose that $\Omega^s$ is a $2$-form on $\k^*$ satisfying the condition~\eqref{eq: Omega s def}. Let $\mathcal{D}$ be the delinearization operator with respect to this choice of $\Omega^s$ and the given $(K,0)$-action on $M$.  Let $(M,\omega) \mod_0 H$ and $(M,\mathcal{D}(\omega)) \mod_0 H$ denote the symplectic quotients of $(M,\omega)$ and $(M,\mathcal{D}(\omega))$, respectively, by the  Hamiltonian $(H,0)$-action with respect to the moment map $\nu$. 
	Then 
	\[
	\widetilde{\mathcal{D}(\omega)} = \mathcal{D}(\widetilde{\omega})
	\]
	as $2$-forms on the underlying manifold $\nu^{-1}(0)/H$ of the symplectic quotients. 
	In particular, the action of $K$ on $M$ descends to a Hamiltonian action of $(K,\pi_{K,s})$ on $(M, \mathcal{D}(\omega))\mod_0 H$ with moment map $\widetilde \Psi$ defined by the property that 
	\[
		\widetilde \Psi \circ pr = \Psi \circ \iota.
	\]
\end{proposition}

\begin{proof} 
	The quotient map $pr\colon \nu^{-1}(0) \to \nu^{-1}(0)/H$ is a surjective submersion. Hence, it suffices to show that 
	\begin{equation}\label{eq: delin quotients} 
	pr^* \left( \widetilde{\mathcal{D}(\omega)}\right) = pr^* \left( \mathcal{D}(\widetilde{\omega}) \right). 
	\end{equation} 
	By definition, the left side of~\eqref{eq: delin quotients} is $pr^*  \left( \widetilde{\mathcal{D}(\omega)}\right) = \iota^* (\mathcal{D}(\omega))$
	where $\iota: \nu^{-1}(0) \into M$ is the inclusion map mentioned above. By definition of delinearization we obtain 
	\[
	\iota^* (\mathcal{D}(\omega)) = \iota^*(\omega + \mu^* \Omega^s) = \iota^* \omega + \iota^* \mu^* \Omega^s = \iota^{*}\omega + (\mu \circ \iota)^{*}\Omega^s.
	\]
	Looking now at the right side of~\eqref{eq: delin quotients}, we may compute 
		\[
	 pr^* \left( \mathcal{D}(\widetilde{\omega}) \right) = pr^*(\widetilde{\omega} + \widetilde{\mu}^* \Omega^s) = pr^*(\widetilde{\omega}) + pr^*(\widetilde{\mu}^* \Omega^s) = \iota^* \omega + (\widetilde{\mu}\circ pr)^* \Omega^s
	 \]
	 where we have used the definition of delinearization in the first equality. By definition of $\widetilde{\mu}$, we have $\mu \circ \iota = \widetilde{\mu} \circ pr$, which completes the proof that $\widetilde{\mathcal{D}(\omega)} = \mathcal{D}(\widetilde{\omega})$.
	 
	 There is a Hamiltonian $(K,0)$-action on $(M,\omega) \mod_0 H$ with moment map $\widetilde \mu$ defined by the property that $\widetilde \mu \circ pr = \mu \circ \iota$. 
	 Applying Theorem~\ref{theorem: delinearization}, it follows that the induced action of $K$ defines a Hamiltonian action of $(K,\pi_{K,s})$ with respect to the symplectic structure $\widetilde{\mathcal{D}(\omega)} = \mathcal{D}(\widetilde{\omega})$ on $(M,\mathcal{D}(\omega))\mod_0 H$ with moment map $\widetilde \Psi = E_s \circ \widetilde \mu$. Since $\widetilde{\mu}$ is defined by the property $\widetilde{\mu} \circ pr = \mu \circ \iota$, it follows that $\widetilde{\Psi}$ is defined by the property $\widetilde \Psi \circ pr = \Psi \circ \iota$ 
	 where $\Psi = E_s \circ \mu$. This completes the proof. 
	 \end{proof}

\section{A local normal form for Poisson-Lie group actions with $K^*$-valued moment maps}\label{sec: local}

The main result of this section is a ``local normal form'' theorem for Hamiltonian actions of the semisimple Poisson-Lie groups $(K,\pi_{K,s})$. This is in the spirit of the Marle-Guillemin-Sternberg local normal form for Hamiltonian actions of $(K,0)$, as recounted in Theorem~\ref{theorem: local normal form}. In fact, our local normal form is derived from Theorem~\ref{theorem: local normal form} together with our ``delinearization commutes with quotient'' result (Proposition~\ref{prop: delinearization quotient}) of Section~\ref{sec: commutes}.  We explain this in detail below. 

Throughout this section, $K$ is a compact connected semisimple Lie group,  $s$ is a positive real constant, and $\pi_{K,s}$ denotes the Poisson-Lie group structure defined in Section \ref{subsec: PL background}. We fix a choice of $2$-form $\Omega^s$ satisfying~\eqref{eq: Omega s def}. Let $(M,\omega)$ be a symplectic manifold equipped with a Hamiltonian action of $(K,\pi_{K,s})$. For technical simplicity, we consider local normal forms near points $p \in M$ such that $\Psi(p)$ is the identity $e \in K^*$.  This is analogous to the simplifying assumption we made in Section~\ref{subsec: local normal forms background} that $\mu(p)=0$.

The following is straightforward from the definitions. 

\begin{lemma} 
In the setting above, we have $T_p(K \cdot p)^{\omega} = \ker(d\Psi)$. 
\end{lemma}

It is a classical fact that $T_p(K \cdot p)^\omega = \ker d\mu$ also for moment maps of Hamiltonian $(K,0)$-actions.  Here we let $\mu$ denote the moment map for the Hamiltonian action of $(K,0)$ on the linearized symplectic manifold $(M,\mathcal{L}(\omega))$, defined with respect to the fixed choice of $2$-form $\Omega^s$ above. Since $\Psi = E_s \circ \mu$ by construction and $E_s$ is a diffeomorphism, we have $\ker(d \Psi) = \ker(d\mu)$.  Thus,
\begin{equation}\label{eq; kernels}
	T_p(K \cdot p)^{\omega} = \ker(d\Psi) = \ker(d\mu) = T_p(K \cdot p)^{\mathcal{L}(\omega)}.
\end{equation}
It follows that 
\begin{equation}\label{eq: symplectic cross sections equal}
	T_p(K\cdot p)^{\mathcal{L}(\omega)} /( T_p(K\cdot p)^{\mathcal{L}(\omega)} \cap  T_p(K\cdot p))= T_p(K\cdot p)^\omega /( T_p(K\cdot p)^\omega \cap  T_p(K\cdot p)).
\end{equation} 
We let $W$ denote the quotient vector space in~\eqref{eq: symplectic cross sections equal} and call it the symplectic cross-section at $p$. Next we consider the two symplectic structures that can be defined on $W$:  one induced by $\omega$, the other induced by $\mathcal{L}(\omega)$. It is useful to know that these symplectic structures coincide, as we show in the next lemma.

\begin{lemma} \label{lem; symplectic cross section} 
In the setting above, $\omega \vert_{T_p(K \cdot p)} = \mathcal{L}(\omega) \vert_{T_p(K \cdot p)}$.
\end{lemma} 

\begin{proof} 
By the definition of $\mathcal{L}(\omega)$ and \eqref{eq; kernels},
\begin{align*}
\omega\vert_{T_p(K\cdot p)^{\omega}}&=(\mathcal{L}(\omega) + \mu^{*}\Omega^s)\vert_{T_p(K\cdot p)^{\omega}}
\\
& = \mathcal{L}(\omega)\vert_{T_p(K\cdot p)^{\omega}} + \Omega^s(d\mu(-)) \vert_{T_p(K \cdot p)^{\omega}} \\
&= \mathcal{L}(\omega)\vert_{T_p(K\cdot p)^{\omega}}
\end{align*}
where in the last equality we have used that $T_p(K \cdot p)^{\omega} = \ker d\mu$. 
From the definition~\eqref{eq: symplectic cross sections equal} of $W$ it follows that $\omega$ and $\mathcal{L}(\omega)$ 
induce the same structure on the quotient vector space $W$. 
\end{proof} 

By Lemma~\ref{lem; symplectic cross section} we may unambiguously denote the symplectic structure on $W$ induced by either $\omega$ or $\mathcal{L}(\omega)$ by $\omega_W$. Now we recall from Section \ref{subsec: local normal forms background} that $(T^*K,\omega_{can})$ is equipped with a Hamiltonian action of $(K\times K,0)$ induced from left and right multiplication on $K$, with moment map $(\mu_L, \mu_R)$ as in~\eqref{eq: muL and muR}. 
Let $(T^*K,\mathcal{D}(\omega_{can}))$ denote the delinearization of $(T^*K,\omega_{can})$ with respect to the moment map $\mu_L$ amd the fixed $2$-form $\Omega^s$. The symplectic manifold $(T^*K,\mathcal{D}(\omega_{can}))$ is equipped with a Hamiltonian action of $(K,\pi_{K,s})$ with moment map $\Psi_L = E_s \circ \mu_L$. By Lemma \ref{lemma: delinearized Hamiltonian by H}, the right action of $(K,0)$ on $(T^*K, \mathcal{D}(\omega_{can}))$ is Hamiltonian, with moment map $\mu_R$.

We can now state our local normal form theorem. Let $H := K_p$ be the stabilizer of $p$. Recall that the symplectic cross-section $W$ is equipped with the isotropy representation of $H$. We may form the product symplectic manifold 
\[
(T^*K \times W, \mathcal{D}(\omega_{can}) \oplus \omega_W)
\]
equipped with an $(H,0)$-action given by $h \cdot (k,\lambda,w) = (kh^{-1},Ad_h^*\lambda,h\cdot w)$. This is a Hamiltonian action of $(H,0)$ with moment map $\mu_R + \Phi_W$.  Our local model is then the symplectic quotient 
\begin{equation}
	(T^*K \times W, \mathcal{D}(\omega_{can}) \oplus \omega_W)\mod_0 H.
\end{equation}
Since the $(K,0)$-action and $(H=K_p,0)$-action arise from a product $(K \times K, 0)$-action on $T^*K \times W$, it follows by Proposition \ref{prop: delinearization quotient} that the action of $(K, \pi_{K,s})$ on $(T^*K \times W, \mathcal{D}(\omega_{can}) \oplus \omega_W)$ descends to a Hamiltonian action of $(K,\pi_{K,s})$ with moment map
\[
	\widetilde \Psi_L \colon (T^*K \times W, \mathcal{D}(\omega_{can}) \oplus \omega_W)\mod_0 H \to K^*, \quad \widetilde \Psi_L \circ pr = \iota \circ \Psi_L
\]
where $\iota$ is the inclusion map for the 0 level set of the $H$ moment map. 

Our main result in this section is the following. 

\begin{theorem}\label{theorem: local normal form for Lu mom maps} 
Let  $K$ be a compact connected semisimple Lie group and $(M,\omega)$ be a symplectic manifold. Suppose that $(M,\omega)$ is equipped with a Hamiltonian $(K, \pi_{K,s})$-action with moment map $\Psi: M \to K^*$. Let $s>0$ be a positive real parameter and suppose that $\Omega^s$ is a $2$-form on $\k^*$ satisfying the condition~\eqref{eq: Omega s def}. Let $\mathcal{D}$ be the delinearization operator with respect to this choice of $\Omega^s$. Then 
there exists a $K$-equivariant symplectomorphism $\varphi$ from a $K$-invariant neighbourhood of  $[e,0,0] $ in $(T^*K \times W, \mathcal{D}(\omega_{can}) \oplus \omega_W)\mod_0 H$ to a $K$-invariant neighbourhood of $p$ in $(M,\omega)$ such that $\varphi([e,0,0]) = p$ and $\Psi \circ \varphi = \widetilde \Psi_L$.
\end{theorem}

Before embarking on the proof, we record a diagram below, which may help the reader visualize the relationships among the different constructions being considered. 
\begin{equation}\label{eq: comm diag} 
\xymatrix @C=5pc { 
(T^*K \times W, \omega_{can} \oplus \omega_W) \ar[r]^{\mathcal{D} \textup{ w.r.t. } K} \ar[d]_{\textup{ quotient at $0$ by $H=K_p$}} & 
(T^*K \times W, \mathcal{D}(\omega_{can}) \oplus \omega_W) \ar[d]^{\textup{ quotient at $0$ by $H=K_p$}} \\ 
(T^*K \times W, \omega_{can} \oplus \omega_W) \mod_0 K_p \ar[r]_{\mathcal{D} \textup{ w.r.t. } K} & 
(T^*K \times W, \mathcal{D}(\omega_{can}) \oplus \omega_W) \mod_0 K_p \\
}
\end{equation}

Roughly, the idea of Theorem~\ref{theorem: local normal form for Lu mom maps} is as follows. We already know from the classical local normal form, Theorem~\ref{theorem: local normal form}, that the lower-left corner of the diagram~\eqref{eq: comm diag} gives a local normal form for the Hamiltonian $(K,0)$-action of the \emph{linearized} symplectic structure $\mathcal{L}(\omega)$. The bottom rightward arrow in the diagram~\eqref{eq: comm diag}, which delinearizes this local normal form, therefore gives us a local normal form for the $(K,\pi_{K,s})$-action with the original moment map $\Psi$ (since $\mathcal{D} \circ \mathcal{L}(\omega) = \omega$ by construction). On the other hand, by Proposition~\ref{prop: delinearization quotient} we know that the diagram~\eqref{eq: comm diag} is commutative. Thus, another equally valid description of this local normal form is given by taking the symplectic quotient of $(T^*K \times W, \mathcal{D}(\omega) \oplus \omega_W)$ by $H=K_p$. This is the content of Theorem~\ref{theorem: local normal form for Lu mom maps}.

\begin{proof} 
Applying linearization to the Hamiltonian action of $(K,\pi_{K,s})$ on $(M,\omega)$ produces a symplectic manifold  $(M,\mathcal{L}(\omega))$ equipped with a Hamiltonian action of $(K,0)$ with moment map $\mu = E_s^{-1} \circ \Psi$.  

The model for the Marle-Guillemin-Sternberg local normal form at a point $p\in M$ for the Hamiltonian action of $(K,0)$ on $(M, \mathcal{L}(\omega)$ is the symplectic quotient  
\[
	(T^*K \times W,\omega_{can}\oplus \omega_W)\mod_0 K_p 
\]
By Theorem \ref{theorem: local normal form}, there exists a $K$-equivariant map $\varphi$ from a $K$-invariant neighbourhood of  $[e,0,0] \in (T^*K \times W, \omega_{can} \oplus \omega_W) \mod_0 K_p$ to a $K$-invariant neighbourhood of $p \in M$ such that $\varphi([e,0,0]) = p$ and $\mu\circ \varphi = \widetilde\mu_L$. Moreover, the map $\varphi$ is a symplectic isomorphism with respect to $\mathcal{L}(\omega)$, i.e. 
\[
	\varphi^*\mathcal{L}(\omega) = \widetilde{\omega_{can}\oplus \omega_W}.
\]
It follows immediately that 
\[
	\Psi \circ \varphi = E_s \circ \mu \circ \varphi = E_s \circ \widetilde\mu_L.
\]
One sees from the defining properties that $\widetilde \Psi_L = E_s \circ \widetilde\mu_L$.

It remains to show that $\varphi$ is a symplectic isomorphism with respect to the original symplectic structure $\omega$ on $M$ and the symplectic structure $\widetilde{\mathcal{D}(\omega_{can}) \oplus \omega_W}$ constructed above. Indeed,
\begin{equation}
\begin{split}
	\varphi^*\omega &  = \varphi^*\mathcal{D}(\mathcal{L}(\omega)) \\
	& = \varphi^*(\mathcal{L}(\omega) + \mu^*\Omega^s )\\
	& = \varphi^*\mathcal{L}(\omega) + \varphi^*\mu^*\Omega^s \\
	& = \widetilde{\omega_{can}\oplus \omega_W} + \widetilde\mu_L^*\Omega^s
\end{split}
\end{equation}
 The last expression is precisely the delinearization of the reduced symplectic structure
 \[
 	\mathcal{D}\left(\widetilde{\omega_{can}\oplus \omega_W}\right) =  \widetilde{\omega_{can}\oplus \omega_W} + \widetilde\mu_L^*\Omega^s
 \]
with respect to the action of $K$. Since delinearization commutes with reduction (Proposition \ref{prop: delinearization quotient}),
 \[
 	\mathcal{D}\left(\widetilde{\omega_{can}\oplus \omega_W}\right) 
	=  \widetilde{\mathcal{D}(\omega_{can}\oplus \omega_W)}. 
 \]
Finally, we have that 
\[
	\mathcal{D}(\omega_{can}\oplus \omega_W) = \mathcal{D}(\omega_{can})\oplus \omega_W
\]
since $d\mu_L$ vanishes on $TW \subset T(T^*K \times W)$. Combining this with the previous identity completes the proof.
\end{proof} 

In summary, a local description of a symplectic manifold $(M,\omega)$ equipped with a Hamiltonian action of $(K,\pi_{K,s})$ with moment map $\Psi$ in a neighbourhood of a point $p$ such that $\Psi(p) = e$ can be obtained from the following data.
\begin{enumerate}
	\item The isotropy subgroup $H = K_p$.
	\item The symplectic cross-section $W = T_p(K\cdot p)^\omega /( T_p(K\cdot p)^\omega \cap  T_p(K\cdot p))$ along with the induced symplectic structure $\omega_W$ and the isotropy representation of $H$.
\end{enumerate}
In particular, the symplectic cross-section  $W$ and the symplectic structure $\omega_W$ can be computed either  from the original symplectic structure $\omega$ or its linearization $\mathcal{L}(\omega)$. As we have observed above, the result is the same. Thus, in order to understand the symplectic structure $\omega$, all that remains is to understand the delinearization of the canonical symplectic structure $\mathcal{D}(\omega_{can})$ on the cotangent bundle $T^*K$. We address this problem in Sections~\ref{sec: delin formula} and~\ref{sec: SU2}.

\section{The symplectic Poisson structure of $\mathcal{D}(\omega_{can})$}\label{sec: delin formula}

In the previous section, we derived a local normal form for Poisson-Lie group actions with $AN$-valued moment maps. The delinearization, $\mathcal{D}(\omega_{can})$, of the canonical symplectic structure on $T^*K$ was an important ingredient in this result.  This leads naturally to the problem of computing a more explicit expression for $\mathcal{D}(\omega_{can})$. The goal of this section and Section~\ref{sec: SU2} is to prove some results in this direction. 

We first sketch the overall strategy and results of this section. Recall that any symplectic structure can equivalently be encoded by its symplectic Poisson structure.  
In this section, we use this perspective and describe the symplectic Poisson structure of  $\mathcal{D}(\omega_{can})$ which we denote by 
\begin{equation}\label{eq: pi def} 
\pi \in \Gamma(\wedge^2 T(T^*K)). 
\end{equation} 
As with $\mathcal{D}(\omega_{can})$, this structure is $s$-dependent, but we suppress $s$ from the notation.
Our main results in this section, Propositions~\ref{prop; symp poisson str coefficients} and~\ref{prop; symplectic poisson structure final}, compute an expression for the coefficients of $\pi$ with respect to the frame 
of $\wedge^2 T(T^*K)$ dual to the frame \eqref{eq; wedge 2 cotangent frame}  which we used to express $\omega_{can}$ in Section~\ref{subsec: local normal forms background}.  Our approach uses the moment map equation for the left and right actions of $K$. Corresponding to the fact that the delinearization $\mathcal{D}(\omega_{can})$ is not unique, we are not able to compute all coefficients of $\pi$. An exception is the special case $K=SU(2)$, where the delinearization is unique and we are able to solve for all coefficients. We treat this special case in  Section~\ref{subsec: pi for SU2}.

\subsection{A formula for $\pi$ at points in the cross-section}\label{subsec: cross section} 

We begin our analysis by examining $\pi$ at points of the form $(e, \lambda) \in K \times \k^*$ where $e \in K$ is the group identity and $\lambda$ is an element of the positive Weyl chamber $\t^*_+ \subseteq \k^*$.
These points form a cross-section for the $K \times K$-action on $T^*K$. In Section~\ref{sec: formula at arbitrary points} we will use this observation to obtain results about $\pi$ at arbitrary points. 

We begin by describing the frame in which we will evaluate coefficients of $\pi$. Recall from Section~\ref{subsec: L background} we have fixed a choice of basis 
$\mathcal{B}\subset \k$. We used these bases to define a frame \eqref{eq; wedge 2 cotangent frame} of $\wedge^2 T^*(T^*K)$. The dual frame of $\wedge^2 T(T^*K)$ in which we will express $\pi$ is
\begin{equation}\label{eq; wedge 2 tangent frame}
	\{ a^L\wedge b^L \mid a,b \in \mathcal{B}\} \cup \{ a^L\wedge b^* \mid  a \in \mathcal{B},\, b^* \in \mathcal{B}^*\} \cup \{ a\wedge b \mid  a^*,b^* \in \mathcal{B}^*\}.
\end{equation}
With respect to this frame, $\pi$ can be written as 
\begin{equation}\label{eq pi general form}
	\pi = \frac{1}{2}\sum_{a,b\in \mathcal{B}} \pi_{a,b} a^L \wedge b^L + \sum_{a,b \in \mathcal{B}} \pi_{a,b^*}a^L\wedge b^* + \frac{1}{2}\sum_{a,b\in \mathcal{B}} \pi_{a^*,b^*}a^* \wedge b^*
\end{equation}
where the coefficients are all skew-symmetric in their arguments. Note that the coefficients $\pi_{a,b}$, $\pi_{a,b^*}$, and $\pi_{a^*,b^*}$ are functions on $T^*K$. The first and last sum in this expression have a factor of $1/2$ because we have not ordered the bases $\mathcal{B}$ and $\mathcal{B}^*$ (and wedge product is skew symmetric). Evaluating this expression at a point of the form $(e,\lambda)$ we obtain
\begin{equation}\label{eq pi on slice}
	\pi(e,\lambda) = \frac{1}{2}\sum_{a,b\in \mathcal{B}} \pi_{a,b}(e,\lambda)a \wedge b + \sum_{a,b \in \mathcal{B}} \pi_{a,b^*}(e,\lambda)a\wedge b^* + \frac{1}{2}\sum_{a,b\in \mathcal{B}} \pi_{a^*,b^*}(e,\lambda)a^* \wedge b^*.
\end{equation}

Our first result computes the coefficients contained in the middle term and the last term of the RHS of~\eqref{eq pi general form} 
by using the moment map $\mu_R$ for the right action of $K$ on $T^*K$. We have the following.

\begin{lemma}\label{lemma: muR comp} 
Let $(e,\lambda) \in T^*K$ where $\lambda \in \t_+^*$. Let $\pi$ be the Poisson bivector field corresponding to $\mathcal{D}(\omega_{can})$. Then, for any  $a,b \in \mathcal{B}$ we have that
\[
	\pi_{a,b^*}(e,\lambda) = \delta_{a,b}, \quad \pi_{a^*,b^*}(e,\lambda) =  -\langle\lambda,[a,b]\rangle.
\]
\end{lemma}

\begin{proof}
From Section~\ref{subsec: local normal forms background} we know that $K \times K$ acts Hamiltonianly on $(T^*K, \omega_{can})$ with moment map $(\mu_R, \mu_L)$. From Lemma \ref{lemma: delinearized Hamiltonian by H} we know that $\mu_R$ satisfies the moment map equation~\eqref{k* valued moment map} with respect to $\mathcal{D}(\omega_{can})$, the delinearization with respect to $\mu_L$. 
Thus the moment map equation
\[
\mathcal{D}(\omega_{can})(a_{T^*K}, \cdot) = d \langle \mu_R, a \rangle
\]
can be rewritten as 
\[
a_{T^*K} = -\pi^\sharp(d \langle \mu_R, a \rangle)
\]
where $a_{T^*K}$ denotes the infinitesmal generated vector field of $a \in \k$ on $T^*K$ for a choice of element $a \in \mathcal{B}$. By~\eqref{eq: K by K action on T*K} it follows that 
this infinitesmal generated vector field of $a$ at the point $(e,\lambda)$ is 
$(-a, ad_a^* \lambda)$. Therefore,
	\begin{equation}\label{lem 1 moment map eq}
		(a,-ad_a^*\lambda) = \pi^\sharp(d\langle \mu_R,a\rangle)(e,\lambda)
	\end{equation}
	which is valid for all $a\in \mathcal{B}$. 
	
	From~\eqref{eq: muL and muR} it is straightforward to compute that, 
	for arbitrary $(X, \eta) \in T_{(e,\lambda)}(K \times \k^*)$, we have 
	$d\langle \mu_R, a\rangle(X, \eta) = - \langle \eta, a \rangle.$
	 If $\eta$ is an element of the dual basis $\mathcal{B}^*$ then this pairing equals $-1$ if $\eta = a^*$ and 0 otherwise. Thus, evaluating the RHS of \eqref{lem 1 moment map eq} at $(e,\lambda)$ yields
	\begin{equation*}
		\pi^\sharp(d\langle \mu_R,a\rangle)(e,\lambda)  = \pi^\sharp_{(e,\lambda)}(-a) 
		 = \sum_{b \in \mathcal{B}}\pi_{b,a^*}(e,\lambda) b -  \sum_{b\in \mathcal{B}} \pi_{a^*,b^*}(e,\lambda)b^*.
	\end{equation*}
	Plugging this back into \eqref{lem 1 moment map eq} yields the system of equations
	\begin{align}
		a  = \sum_{b \in \mathcal{B}}\pi_{b,a^*}(e,\lambda) b, \qquad -ad_a^*\lambda  = -\sum_{b\in \mathcal{B}} \pi_{a^*,b^*}(e,\lambda)b^*.
	\end{align}
	It follows from the first equation that $\pi_{b,a^*}(e,\lambda) =  \delta_{a,b}$. The second equation tells us that	
	\[
		\pi_{a^*,b^*}(e,\lambda) = \langle ad_a^* \lambda, b\rangle = -\langle \lambda, [a,b] \rangle.  \qedhere
	\]

\end{proof}

From Lemma~\ref{lemma: muR comp} it immediately follows that 
\begin{equation}\label{eq: pi e lambda 2} 
	\pi(e,\lambda) = \frac{1}{2}\sum_{a,b\in \mathcal{B}} \pi_{a,b}(e,\lambda) a \wedge b + \sum_{b \in \mathcal{B}} b\wedge b^* - \frac{1}{2}\sum_{a,b\in \mathcal{B}} \langle\lambda,[a,b]\rangle a^* \wedge b^*.
\end{equation}

We now turn our attention to computing $\pi_{a,b}(e,\lambda)$. To 
analyze these coefficients we 
use the $AN$-valued moment map $\Psi= E_s \circ \mu_L$. By an argument similar to the preceding proof, the moment map equation can be reformulated as 
\[
a_{T^*K} = -\pi^\sharp(\Psi^* \langle \Theta^L, a \rangle_s)
\]
for all $a \in \mathcal{B}$, where $a_{T^*K}$ now denotes  the infinitesmal generated vector field  with respect to the \emph{left} multiplication action of $K$. (In the preceding proof, this  notation was used for the infinitesmal generated vector field with respect to the \emph{right} multiplication action.) By~\eqref{eq: K by K action on T*K}, the value of this vector field at $(e,\lambda)$ is equal to $(a,0)$, so we conclude 
	\begin{equation}\label{lem 2 moment map eq}
		(a,0) = -\pi^\sharp(\Psi^*\langle \Theta^L,a\rangle_s)(e,\lambda)
	\end{equation}
	for all elements $a \in \mathcal{B}$.

	We now analyze the RHS of~\eqref{lem 2 moment map eq}, and in particular,  derive an expression for the covector $(\Psi^*\langle \Theta^L,a\rangle_s)_{(e,\lambda)}$.  The first computation, recorded in Lemma~\ref{lemma: 1-form easy half}, gives the pairing of this covector 
	against vectors of the form $(b,0)$. A similar computation also appears in \cite[Proposition 5.12]{Lu}.

\begin{lemma}\label{lemma: 1-form easy half}
Let $a,b \in \mathcal{B}$ and let $\lambda \in \t_+^*$.  Then,
	\[
		( \Psi^*\langle \Theta^L,a \rangle_s)_{(e,\lambda)}(b,0) = \sum_{\alpha \in R_+} \frac{1 - e^{-2s\sqrt{-1}\alpha(\phi(\lambda))}}{2s}  \left(\delta_{b,y_\alpha}\delta_{a,x_\alpha} - \delta_{b,x_\alpha}\delta_{a,y_\alpha} \right).
	\]
\end{lemma}

\begin{proof}  From~\eqref{eq: K by K action on T*K} we see that, as an element of the tangent space $T_{(e,\lambda)}(K\times \k^*)$,  $(b,0)$ is the infinitesmal generated vector field of the \emph{left} action of $K$ on $T^*K$. It  follows by \eqref{K* valued moment map} and Lemma \ref{background useful fact}  that
	\begin{equation}\label{lem 2 moment map eq 2}
		(\Psi^*\langle \Theta^L,a\rangle_s)_{(e,\lambda)}(b,0) = (L_{p^{-1}})_*(\pi_{AN,s})_p(a^*,b^*)
	\end{equation}
	where 
	\[
		p =\Psi(e,\lambda) = E_s(\lambda) = \exp(s\sqrt{-1}\phi(\lambda)).
	\] 
	Thus, our problem reduces to evaluating the RHS of~\eqref{lem 2 moment map eq 2}. From Lemma \ref{lemma: piAN take 1}
	it follows that  
	\[
		(L_{p^{-1}})_*(\pi_{K^*,s})_p = s\sum_{\alpha \in R_+} (p^{-2\alpha}-1) \sqrt{-1}e_\alpha\wedge e_\alpha.
	\]
	From this it is easy to see that \eqref{lem 2 moment map eq 2} will vanish for all combinations of $a,b\in \mathcal{B}$ except when $a=x_\alpha$ and $b=y_\alpha$ (or vice versa) for some $\alpha \in R_+$. We compute 
	\begin{align*}
		((L_{p^{-1}})_*(\pi_{K^*,s})_p)(x_\alpha^*,y_\alpha^*) & = s(p^{-2\alpha}-1) (\langle \sqrt{-1}e_\alpha,x_\alpha\rangle_s\langle e_\alpha,y_\alpha\rangle_s - \langle \sqrt{-1}e_\alpha,y_\alpha\rangle_s \langle e_\alpha,x_\alpha \rangle_s) \\
		& = \frac{1-p^{-2\alpha}}{2s}.
	\end{align*}
	Finally, for $p = \exp(s\sqrt{-1}\phi(\lambda))$ we have 
	$p^{-2\alpha} = e^{-2\alpha(s\sqrt{-1}\phi(\lambda))} = e^{-2s\sqrt{-1}\alpha(\phi(\lambda))}$.
	The lemma follows.
\end{proof}

It remains to show how the covector $(\Psi^*\langle \Theta^L,a\rangle_s)_{(e,\lambda)}$ evaluates on tangent vectors of the form $(0,b^*)$. We introduce the notation 
\begin{equation}\label{beta coefficients definition}
	\beta_{a,b} := (\Psi^*\langle \Theta^L,a\rangle_s)_{(e,\lambda)}(0,b^*) = \langle \Theta^L(\Psi_*)_{(e,\lambda)}(0,b^*),a\rangle_s = \langle \Theta^L(dE_s)_\lambda(b^*),a\rangle_s
\end{equation}
for all $a,b$ in $\mathcal{B}$.  By the definition of the $\beta_{a,b}$ and Lemma~\ref{lemma: 1-form easy half}, at a point of the form $(e,\lambda)$ in $T^*K$ we have 
\begin{equation}\label{eq: beta a b equation} 
	(\Psi^*\langle \Theta^L,a\rangle_s)_{(e,\lambda)} = \sum_{\alpha \in R_+}\frac{1 - e^{-2s\sqrt{-1}\alpha(\phi(\lambda))}}{2s}\left(  \delta_{a,x_\alpha}y_\alpha^*-\delta_{a,y_\alpha}x_\alpha^*\right)+ \sum_{b \in \mathcal{B}}\beta_{a,b}b
\end{equation} 
as an element of the cotangent space $T^*_{(e,\lambda)}(K \times k^*) \cong \k^* \times k$. The following lemma gives a collection of linear equations which must be satisfied by the coefficients $\{\beta_{a,b}\}$. 

\begin{lemma}\label{lemma: system of eqs for coefficients} 
Let the notation be as above. Let $a \in \mathcal{B}$.  Then 
\begin{align}\label{eq: system of eqs for coefficients}
	a & =  \sum_{b \in \mathcal{B}} \beta_{a,b}b + \left(\frac{1}{2}\sum_{b,c\in \mathcal{B}} \pi_{b,c}(e,\lambda) b \wedge c  \right)^\sharp\sum_{\alpha \in R_+}\frac{1 - e^{-2s\sqrt{-1}\alpha(\phi(\lambda))}}{2s}\left(  \delta_{a,x_\alpha}y_\alpha^*-\delta_{a,y_\alpha}x_\alpha^*\right) , \\
	0 & =  \sum_{\alpha \in R_+}\frac{1 - e^{-2s\sqrt{-1}\alpha(\phi(\lambda))}}{2s}\left(  \delta_{a,x_\alpha}y_\alpha^*-\delta_{a,y_\alpha}x_\alpha^*\right) +\left( \frac{1}{2}\sum_{b,c\in \mathcal{B}}\langle \lambda,[b,c]\rangle b^* \wedge c^*\right)^\sharp \sum_{b \in \mathcal{B}}\beta_{a,b}b. 
\end{align}
\end{lemma} 

\begin{proof} 
It suffices to plug the expression~\eqref{eq: beta a b equation} into the moment map equation~\eqref{lem 2 moment map eq}. 
\end{proof}

Summarizing the results obtained thus far, we have the following. 

\begin{proposition}\label{prop; symp poisson str coefficients}
Let $(e,\lambda) \in T^*K$ where $\lambda \in \t_+^*$. Let $\pi$ be the Poisson bivector field corresponding to $\mathcal{D}(\omega_{can})$.
Then \[
\pi(e,\lambda) = \frac{1}{2}\sum_{a,b\in \mathcal{B}} \pi_{a,b}(e,\lambda) a \wedge b + \sum_{b \in \mathcal{B}} b\wedge b^* - \frac{1}{2}\sum_{a,b\in \mathcal{B}} \langle\lambda,[a,b]\rangle a^* \wedge b^*.
\]
Moreover, the coefficients $\pi_{a,b}(e,\lambda)$ must satisfy the linear equations~\eqref{eq: system of eqs for coefficients}, where the $\beta_{a,b}$ are defined by~\eqref{beta coefficients definition}. 
\end{proposition}

In the general situation, the equations~\eqref{eq: system of eqs for coefficients} does not fully specify the coefficients $\pi_{a,b}(e,\lambda)$. As we will see in Section~\ref{subsec: pi for SU2}, these equations  completely determine the symplectic Poisson structure in the special case $K=SU(2)$.  

\subsection{A formula for $\pi$ at arbitrary points in $T^*K$} \label{sec: formula at arbitrary points}

In the previous section, we described the coefficients of $\pi$ at points of the form $(e,\lambda) \in T^*K$. In this section, we use this information to described the coefficients of  $\pi$ at arbitrary points in $T^*K$. To do this, we use the 
 fact that the left multiplication action of $K$ on $(T^*K ,\pi)$ is Poisson with respect to $\pi_{K,s}$ and the right multiplication action of $K$ is Poisson with respect to the trivial Poisson structure on $K$ (cf. Section~\ref{subsec: HPL background} and in particular, equation~\eqref{action map is Poisson}). We start with the following. 
\begin{proposition}\label{prop; symplectic poisson structure final}
Let $(k, \xi) \in K \times \k^* \cong T^*K$. The
symplectic Poisson structure $\pi$ evaluates at $(k, \xi)$ to be
	\begin{align}\label{lem 3 eq}
		\pi(k,\xi) & = \frac{1}{2}\sum_{a,b\in \mathcal{B}}\pi_{a,b}(k,\xi)a^L \wedge b^L + \sum_{b \in \mathcal{B}} b^L\wedge  b^* - \frac{1}{2} \sum_{a,b\in \mathcal{B}}  \langle \xi, [a,b] \rangle a^* \wedge b^* 
		+ ( Ad_{kh}r_{K,s} - r_{K,s})^R_k. 
	\end{align}
	where $h \in K$ and $\lambda \in \t^*_+$ are such that $\xi = Ad_{h}^*\lambda$. The coefficients $\pi_{a,b}(k,\xi)$ are related to the coefficients $\pi_{a,b}(e,\lambda)$ by the equation
	\[
		 \frac{1}{2}\sum_{a,b\in \mathcal{B}}\pi_{a,b}(k,\xi)a^L \wedge b^L = \frac{1}{2}\sum_{a,b\in \mathcal{B}} \pi_{a,b}(e,\lambda) (Ad_{h}a)^L \wedge (Ad_{h}b)^L.
	\]
\end{proposition}

\begin{proof}
By delinearization (Theorem~\ref{theorem: delinearization}), we have a Hamiltonian $(K, \pi_{K,s})$-action on $(T^*K, \mathcal{D}(\omega_{can}))$ with moment map $E_s \circ \mu_L$. Throughout we use the identification $T^*K \cong K\times \k^*$ by left invariant covector fields. The action map 
	\begin{equation}\label{eq:  sigma1} 
		\sigma_1 \colon K  \times (K \times \k^*) \to K \times \k^*, \quad \sigma_1(k_1,(k,\xi)) = (k_1k,\xi)
	\end{equation} 
	for the left multiplication action of $K$  is Poisson with respect to the product of the Poisson structure $\pi_{K,s}$ on $K$ and the Poisson structure $\pi$. 
		On the other hand, from Lemma~\ref{lemma: delinearized Hamiltonian by H}, we also know that the action map
	\begin{equation}\label{eq: sigma2} 
		\sigma_2 \colon K  \times (K \times \k^*) \to K \times \k^*, \quad \sigma_2(k_2,(k,\xi)) = (kk_2^{-1},Ad_{k_2}^*\xi) 
	\end{equation} 
	corresponding to the \emph{right} multiplication action of $K$ is Poisson with respect to the product of the trivial Poisson structure on $K$ and the symplectic Poisson structure $\pi$.   
	We now use $\sigma_1$ and $\sigma_2$ to translate 
our previous computations at $(e,\lambda)$ to arbitrary points in $T^*K \cong K \times \k^*$.

	Let $(k,\xi)$ be an arbitrary point in $K \times \k^*$. Let $h \in K$ and $\lambda \in \t^*$ such that $\xi = Ad_{h}^* \lambda$. Let $k_1 := kh$. Then it follows from~\eqref{eq: sigma1} and~\eqref{eq: sigma2} that 
$(k,\xi) = \sigma_1(k_1,\sigma_2(h,(e,\lambda)))$. 
	Since $\sigma_1$ and $\sigma_2$ are Poisson maps, we have by applying \eqref{action map is Poisson} twice 
	 that
	\begin{equation}\label{eq: pi k lambda}
	\begin{split} 
		\pi(k,\lambda) & = \pi(\sigma_1(k_1,\sigma_2(h,(e,\lambda))))\\
		& = \sigma_1(k_1)_*( \sigma_2(h)_*\pi(e,\lambda) +\sigma_2(e,\lambda)_* (0) ) + \sigma_1(h^{-1},\xi)_* \pi_{K,s}(k_1)\\
		& = \sigma_1(k_1)_*( \sigma_2(h)_*\pi(e,\lambda) ) + \sigma_1(h^{-1},\xi)_* \pi_{K,s}(k_1).
	\end{split} 
	\end{equation}
	and we use the facts that $\sigma_2(h, (e, \lambda)) = (h^{-1}, \xi)$ and that the right multiplication action is Poisson with respect to the $0$ Poisson structure on $K$, while the left multiplication action is Poisson with respect to the $\pi_{K,s}$ Poisson structure on $K$.

	We compute the second term on the RHS of~\eqref{eq: pi k lambda} first. Recall from Section~\ref{subsec: PL background} that
	\[
		\pi_{K,s}(k_1) = (r_{K,s}^L - r_{K,s}^R)(k_1) = (L_{k_1})_*r_{K,s}- (R_{k_1^{-1}})_*r_{K,s} = (R_{k_1^{-1}})_*( Ad_{k_1}r_{K,s} - r_{K,s}).
	\]
	We need to compute the pushforward of this bivector along the map $\sigma_1(h^{-1}, \xi)_*$, where $\sigma_1(h^{-1},\xi)$ is the map $K \to K \times k^*, k \mapsto \sigma_1(k, (h^{-1}, \xi)) = (kh^{-1}, \xi)$ and the differential of $\sigma_1(h^{-1},\xi)$ must be computed at the point $k_1 \in K$. From the explicit description of $\sigma_1(h^{-1},\xi)$ we can compute this to be 
	\[
		T_{k_1}K \to T_{k_1h^{-1}}K \times \k, \quad X \mapsto ((R_{h})_*X,0)
	\]
	The map $\sigma_1(h^{-1},\xi)_*\colon \wedge^2T_{k_1}K \to \wedge^2 T_{(k_1h^{-1},\xi)}(K\times \k)$ is the second wedge power of this map. Therefore, the pushforward of the bivector $\pi_{K,s}(k_1) \in \wedge^2 T_{k_1}K$  is 
	\[
	\begin{split} 
		\sigma_1(h^{-1},\xi)_* \pi_{K,s}(k_1) & = ((R_{h})_* (R_{k_1^{-1}})_*( Ad_{k_1}r_{K,s} - r_{K,s}), 0) \\
		& = ((R_{h k_1^{-1}})_*( Ad_{k_1}r_{K,s} - r_{K,s}), 0) \\
		& = ((R_{k^{-1}})_* ( Ad_{k_1}r_{K,s} - r_{K,s}), 0) \\
		& = ( Ad_{k_1}r_{K,s} - r_{K,s})^R_{k} \\ 
		& =( Ad_{kh}r_{K,s} - r_{K,s})^R_k
	\end{split} 
	\]
	where recall $k_1 := kh$. 
	
	Next, we compute the first term in the RHS of~\eqref{eq: pi k lambda}.  First, the map $\sigma_2(h): K \times \k^* \to K \times \k^*$ is given by $(k, \xi) \mapsto (k h^{-1}, Ad^*_{h}\xi)$. From this it follows that $\sigma_2(h)_*: \k = T_eK \to T_{h^{-1}}K$ is given by 
	$a \mapsto (R_{h})_*(a) = (L_{h^{-1}})_* (Ad_{h}(a))$.   From this we compute 
	\begin{align*}
		 &\sigma_1(k_1)_*\sigma_2(h)_*\pi(e,\lambda)\\
		 & = \sigma_1(k_1)_*\sigma_2(h)_* \left(\frac{1}{2}\sum_{a,b\in \mathcal{B}}\pi_{a,b}(e,\lambda) a \wedge b +\sum_{b\in \mathcal{B}}b\wedge b^* - \frac{1}{2}\sum_{a,b\in \mathcal{B}} \langle \lambda, [a,b] \rangle a^* \wedge b^* \right) \, \, \\
		& = \sigma_1(k_1)_* \left( \frac{1}{2}\sum_{a,b\in \mathcal{B}} \pi_{a,b}(e,\lambda) (Ad_{h}a)^L \wedge (Ad_{h}b)^L +\sum_{b\in \mathcal{B}}(Ad_{h}b)^L\wedge Ad_{h}^* b^* - \frac{1}{2} \sum_{a,b\in \mathcal{B}} \langle \lambda, [a,b] \rangle Ad_{h}^*a^* \wedge Ad_{h}^*b^* \right) \\
		& =  \frac{1}{2} \sum_{a,b\in \mathcal{B}}\pi_{a,b}(k,\xi) a^L \wedge b^L +\sum_{b\in \mathcal{B}}(Ad_{h}b)^L\wedge Ad_{h}^* b^* - \frac{1}{2} \sum_{a,b\in \mathcal{B}} \langle Ad_{h^{-1}}^*\xi, [a,b] \rangle Ad_{h}^*a^* \wedge Ad_{h}^*b^*\\
		& =  \frac{1}{2} \sum_{a,b\in \mathcal{B}}\pi_{a,b}(k,\xi) a^L \wedge b^L +\sum_{b\in \mathcal{B}}b^L\wedge b^* - \frac{1}{2} \sum_{a,b\in \mathcal{B}} \langle \xi, [a,b] \rangle a^* \wedge b^*. 
	\end{align*}
	In the penultimate line, we used $\xi = Ad_h^*\lambda$. In the last line we used that $\sum_{b\in \mathcal{B}}(Ad_{h}b)^L\wedge Ad_{h}^* b^* = \sum_{b\in \mathcal{B}}b^L\wedge b^*$ and
	\[
		 \frac{1}{2} \sum_{a,b\in \mathcal{B}} \langle Ad_{h^{-1}}^*\xi, [a,b] \rangle Ad_{h}^*a^* \wedge Ad_{h}^*b^* =  \frac{1}{2} \sum_{a,b\in \mathcal{B}} \langle \xi, [Ad_{h}a,Ad_{h}b] \rangle (Ad_{h}a)^* \wedge (Ad_{h}b)^* = \frac{1}{2} \sum_{a,b\in \mathcal{B}} \langle \xi, [a,b] \rangle a^* \wedge b^*.
	\]
	 Putting the first term together with the second term finishes the proof.
\end{proof}

As we see in the next section, we can give an explicit formula for the coefficients $\pi_{a,b}$ in the case $K=SU(2)$. Although we cannot give a formula for the coefficients $\pi_{a,b}$ in  the general case, we can describe their behaviour in the limit $s\to 0$. Note that, by the linearization theorem, the limit as $s \to 0$ of this Poisson structure recovers the symplectic Poisson structure associated to the canonical symplectic structure $\omega_{can}$.  Comparing with the formula \eqref{omega can explicit} for the symplectic form, we have that
	\begin{align}
		\lim_{s\to 0}\pi(k,\xi) & =  \sum_{b\in \mathcal{B}}b^L\wedge  b^* - \frac{1}{2}  \sum_{a,b\in \mathcal{B}}\langle \xi,[a,b]\rangle a^* \wedge b^*.
	\end{align}
We also know that  $( Ad_{kh}r_{K,s} - r_{K,s})^R_k \to 0$ as $s\to 0$ since $r_{K,s} \to 0$ as $s \to 0$. It follows  that  $\pi_{a,b} \to 0$ as $s\to 0$ for all $a,b\in \mathcal{B}$.

\subsection{Example: $K = SU(2)$}\label{subsec: pi for SU2}

In the previous sections, we derived some general results that yield information about the coefficients of the symplectic Poisson structure $\pi$ associated to the delinearized symplectic structure $\mathcal{D}(\omega_{can})$. In the special case $K=SU(2)$, it turns out that the information obtained in the previous sections is enough to determine $\pi$ uniquely. We record the details below.

For $K=SU(2)$ we have $\mathcal{B} = \{x,y,t\}$, where the elements $x,y,t$ are as defined in \eqref{sl2 orthonormal basis}. 
We denote by $\alpha$ the unique positive simple root. For simplicity in the computations that follow we introduce the notation  
\[
	\gamma := \frac{1 - e^{-2s\sqrt{-1}\alpha(\phi(\lambda))}}{2s}.
\] 
Here, $\lambda$ is an element of the positive Weyl chamber of $\mathfrak{su}(2)^*$, so it equals $\xi t^*$ for some $\xi > 0$. Thus 
\[
	\gamma =  \frac{1 - e^{-2s\sqrt{-1}\alpha(\phi(\xi t^*))}}{2s} = \frac{1 - e^{-2s\sqrt{-1}\xi \alpha (t)}}{2s}  = \frac{1 - e^{s \sqrt{2} \xi}}{2s}.
\]

From Lemma~\ref{lemma: system of eqs for coefficients} we have a system of equations satisfied by the coefficients $\pi_{a,b}$ appearing in the expansion of $\pi$ with respect to the basis $\mathcal{B}$. For $K=SU(2)$ these equations become: 
\begin{align}\label{eq1}
	a & =   \beta_{a,t} t + \beta_{a,x} x + \beta_{a,y} y +  \gamma\left(\pi_{t,x}(e,\lambda) t\wedge x + \pi_{x,y}(e,\lambda) x\wedge y + \pi_{y,t}(e,\lambda)y \wedge t  \right)^\sharp\left(\delta_{a,x}y^* - \delta_{a,y}x^* \right) ,
\end{align}
\begin{align}\label{eq2}
	0 & = \gamma (\delta_{a,x}y^*-\delta_{a,y}x^* ) + \frac{1}{\sqrt{2}}\left(\langle\lambda, t\rangle x^*\wedge y^* +\langle\lambda, x\rangle y^*\wedge t^* + \langle\lambda,y\rangle t^* \wedge x^*\right)^\sharp \left(\beta_{a,t}t  + \beta_{a,x}x  + \beta_{a,y}y \right). 
\end{align}

We analyze~\eqref{eq1} and~\eqref{eq2} separately, starting with~\eqref{eq1}. 
Note that these equations must hold for any choice of $a \in \mathcal{B}$; by setting 
$a=t$, $a=x$, and $a=y$ in \eqref{eq1} it is not hard to see that 
\begin{align}\label{linear system 1.1}
	\beta_{x,y}  & = 0,\quad 
	\beta_{y,x}  = 0,\quad
	\beta_{t,t}  = 1,\quad
	\beta_{t,x}  = 0,\quad 
	\beta_{t,y}  = 0,
\end{align}
\begin{align}\label{linear system 1.2}
	1  & = \beta_{x,x} - \gamma \pi_{x,y}(e,\lambda), \quad
	\beta_{x,t}  =  -\gamma \pi_{y,t}(e,\lambda), \quad
	1  = \beta_{y,y} - \gamma \pi_{x,y}(e,\lambda), \quad
	\beta_{y,t}  = -\gamma \pi_{t,x}(e,\lambda).
\end{align}

Next we simplify~\eqref{eq2}. As already noted, $\lambda = \xi t^*$ for some $\xi >0$. This implies that $\langle \lambda,x\rangle = \langle \lambda,y\rangle = 0$, and $\langle \lambda,t\rangle = \xi$ so \eqref{eq2} simplifies to
\begin{align}\label{linear system 2.1}
	0  & = \gamma (\delta_{a,x}y^*-\delta_{a,y}x^* ) + \frac{\xi}{\sqrt{2}} (x^*\wedge y^*)^\sharp \left(\beta_{a,t}t  + \beta_{a,x}x  + \beta_{a,y}y \right) \\
	& =  \gamma (\delta_{a,x}y^*-\delta_{a,y}x^* ) + \frac{\xi}{\sqrt{2}} (  \beta_{a,x}y^* -  \beta_{a,y} x^*).
\end{align}
which in turn yields 
\begin{align}\label{linear system 2.2}
	\gamma \delta_{a,y} = -\frac{\xi}{\sqrt{2}} \beta_{a,y}, \quad \gamma \delta_{a,x} = -\frac{\xi}{\sqrt{2}} \beta_{a,x}
\end{align}
which is valid for all $a \in \mathcal{B}$. 
Using the relations in \eqref{linear system 2.2} to simplify \eqref{linear system 1.2}, we obtain 
\begin{align}\label{linear system 3}
	\pi_{x,y}(e,\lambda) = \frac{\beta_{x,x}-1}{\gamma} = -\frac{1}{\gamma} - \frac{\sqrt{2}}{\xi} \quad \textup{ and } \quad 
	\pi_{y,t}(e,\lambda) = -\frac{1}{\gamma}\beta_{x,t}, \quad  \pi_{t,x} = -\frac{1}{\gamma}\beta_{y,t}.
\end{align}
It remains to compute $\beta_{x,t}$ and $\beta_{y,t}$.   We have the following. 

\begin{lemma} 
$\beta_{x,t} = \beta_{y,t} = 0$. 
\end{lemma}

\begin{proof} 

From the definition \eqref{beta coefficients definition} we have 
\begin{equation}
	\beta_{x,t} = \langle \Theta^L\Psi_*(0,t^*),x\rangle_s = \langle \Theta^L(dE_s)_\lambda(t^*),x\rangle_s \quad \textup{ and } \quad 
	\beta_{y,t} = \langle \Theta^L\Psi_*(0,t^*),y\rangle_s = \langle \Theta^L(dE_s)_\lambda(t^*),y\rangle_s.
\end{equation}
The restriction $E_s\colon \t^* \to A \subset AN$ is $E_s(\lambda) = \exp(s\sqrt{-1}\varphi(\lambda))$. Since $t^*$ (viewed as a tangent vector at $\lambda$) is tangent to $\t^*$, we can use this restriction to compute the pushforward of $t^*$ under $E_s$. We have:
\begin{align}
	(dE_s)_\lambda(t^*) & =  d(\exp_{s\sqrt{-1}\phi(\lambda)})( s\sqrt{-1} t)  =  e^{s\sqrt{-1}\phi(\lambda)}( s\sqrt{-1} t) = E_s(\lambda)( s\sqrt{-1} t)
\end{align}
where by the last expression we mean the left multiplication of $s\sqrt{-1}t \in A$ by $E_s(\lambda)$. 
Thus, this last expression is the value at $E_s(\lambda)$ of the left-invariant vector field on $A$ whose value at the identity is $s\sqrt{-1} t$. From this we conclude 
\begin{align}
	\Theta^L(dE_s)_\lambda(t^*) = s\sqrt{-1} t
\end{align}
and we obtain 
\begin{equation}
	\beta_{x,t} = \langle s\sqrt{-1} t,x\rangle_s = \frac{1}{s} \im \kappa(s\sqrt{-1} t,x) = 0 \quad \textup{ and } \quad 
	\beta_{y,t} = \langle s\sqrt{-1} t,y\rangle_s = \frac{1}{s} \im \kappa(s\sqrt{-1} t,y) = 0
\end{equation}
as desired. 
 \end{proof}

Summarizing the computations thus far, we have shown the following. 

\begin{proposition} 
Let $K = SU(2)$. Let $(k, \xi) \in SU(2) \times \mathfrak{su}(2)^* \cong T^*SU(2)$. Let $\mathcal{B} = \{x, y, t\}$ be the ordered basis of $\mathfrak{su}(2)$ as chosen in~\eqref{sl2 orthonormal basis}. 
Then the 
symplectic Poisson structure $\pi$ associated to $\mathcal{D}(\omega_{can})$ evaluates at $(k, \xi)$ is 
	\begin{equation}
	\begin{split} 
		\pi(k,\xi) & = -\left( \frac{1}{\gamma} + \frac{\sqrt{2}}{\xi} \right) (Ad_{h}x)^L \wedge (Ad_{h}y)^L +\sum_{b \in \mathcal{B}} b^L\wedge  b^*   - \frac{1}{2} \sum_{a,b\in \mathcal{B}}  \langle \xi, [a,b] \rangle a^* \wedge b^* 
		+  ( Ad_{kh}r_{K,s} - r_{K,s})^R_k \\ 
		\end{split} 
	\end{equation}
	where $h \in K$ and $\lambda \in \t^*_+$ are such that $Ad_{h}^*\xi = \lambda$, and $r_{K,s} :=  s \cdot x \wedge y$. 
\end{proposition} 
 
\begin{remark} In this case we can see directly that the coefficient $\pi_{x,y}(e,\lambda)$ goes to zero as $s \to 0$ and thus confirm that 
the delinearized structure approaches the non-delinearized structure as $s \to 0$. First note that 
\[
	\gamma = \frac{1 - e^{s \sqrt{2} \xi}}{2s} = \frac{-\sqrt{2} \xi + O(s^2)}{2s} = - \frac{\xi}{\sqrt{2}} + O(s).
\]
Thus,
\[
	\pi_{x,y}(e,\lambda) = -\frac{1}{\gamma} - \frac{\sqrt{2}}{\xi} = \frac{1}{ \frac{\xi}{\sqrt{2}} + O(s)} - \frac{\sqrt{2}}{\xi}
\]
which goes to 0 as $s \to 0$.
\end{remark}

\section{Example: $\mathcal{D}(\omega_{can})$ in the case of $K=SU(2)$ in standard coordinates}\label{sec: SU2}

The purpose of this section is to derive, using a method different from Section~\ref{sec: delin formula}, explicit formulas for the coefficients defining the delinearized symplectic $2$-form $\mathcal{D}(\omega_{can})$. Unlike our method in Section~\ref{sec: delin formula}, in this section we do the computation directly, with no reference to the Poisson structure $\pi$. The 2-form $\mathcal{D}(\omega_{can})$ and the Poisson structure $\pi$ are related, essentially by a matrix inverse, so it is of course possible to derive one formula from the other.  We record both computations because we felt that the methods are different and both are instructive.

We first describe in broad strokes the overall strategy employed below. Recall that the delinearized structure $\mathcal{D}(\omega_{can})$ is defined as 
$\mathcal{D}(\omega_{can}) := \omega_{can} + \mu_L^* \Omega^s$. The key property~\eqref{eq: Omega s def} 
tells us how to evaluate $\Omega^s$ on an arbitrary vector against a vector tangent to the coadjoint orbits. However, in general, $\Omega^s$ cannot be uniquely determined from this information. 
The elementary observation recorded in the lemma below forms the basis of the computations in this section. 

\begin{lemma}\label{lemma: SU2} 
	Let $K$ be a compact connected semisimple (real) Lie group and let $\k$ and $\k^*$ denote the Lie algebra and dual of the Lie algebra of $K$ respectively. Assume $\k$ is equipped with a $K$-invariant inner product and $\k^*$ is equipped with the induced inner product. Let $\lambda \in \k^*$. If the coadjoint orbit $\mathcal{O}_\lambda := K \cdot \lambda$ through the point $\lambda$ is of codimension $1$, then the equality~\eqref{eq: Omega s def} uniquely determines the $2$-form $\Omega^s$ at the point $\lambda$. 
\end{lemma}

\begin{proof} 
For any $v \in T_\lambda \k^* \cong \k^*$ we may decompose it into directions tangent to the orbit and orthogonal to the orbit, i.e. 
	$v = X_{\k^*} + u$ for some $X \in \k$ and some $u \in (T_\xi \mathcal{O}_\lambda)^{\perp}$. 
	Now for any pair $v, w \in T_\lambda \k^*$ we have $v = X_{\k^*} + u$ and similarly $w = Y_{\k^*} + u'$ for $Y \in \k$ and $u \in (T_\lambda \mathcal{O}_\lambda)^{\perp}$. Then we see that 
	\[
	\Omega^s(v,w) = \Omega^s(X_{\k^*}, Y_{\k^*}) + \Omega^s(X_{\k^*}, u') + \Omega^s(u, Y_{\k^*}) + \Omega^s(u, u')
	\]
	where (by skew-symmetry) the first three terms may be computed using equation~\eqref{eq: Omega s def}. Now note that the last term must be equal to $0$ since by assumption $\mathrm{dim}(T_\lambda \mathcal{O}_\lambda)^{\perp} = 1$ and thus $u$ and $u'$ must be scalar multiples and $\Omega^s$ is skew-symmetric. Hence $\Omega^s(v,w)$ can be completely described by~\eqref{eq: Omega s def}, as claimed. 
\end{proof} 

In the case of $K=SU(2)$, the coadjoint orbits through any non-zero point in $\mathfrak{su}(2)^*$ are $2$-dimensional spheres and hence codimension-$1$. Thus the above lemma applies. 
We will also use in a crucial way the facts that $E_s$ is $K$-equivariant with respect to the coadjoint action of $K$ on $\k^*$ and the dressing action of $K$ on $AN$, and that the inclusion of a dressing orbit into $AN$ is Hamiltonain $(K, \pi_{K,s})$-action with respect to the Poisson structure $\pi_{AN}$ on $AN$. A glance at the definition of $\mathcal{D}(\omega_{can})$, the equation~\eqref{eq: Omega s def}, Lemma~\ref{lemma: SU2}, and the facts just mentioned, reveals that we need several preliminaries: 
\begin{itemize} 
	\item a computation of $E_s$ (Lemma~\ref{lemma: Es coordinates}); 
	\item a computation of the decomposition $T_\lambda \mathcal{O}(\lambda) \oplus (T_\xi \mathcal{O}(\lambda))^{\perp}$ for points $\lambda \in \k^*$ (Lemma~\ref{lemma: projections}); 
	\item a computation of the Poisson bivector $\pi_{AN}$ on $AN$ (Lemma~\ref{lemma: piAN take 2}). 
\end{itemize} 

We prove each of the above lemmas in turn. 
We use the basis $\mathcal{B} := \{x, y, t\}$ for $\mathfrak{su}(2)$ of~\eqref{sl2 orthonormal basis}. As explained in Example~\ref{sl2 orthonormal basis}, this basis is orthonormal with respect to the pairing~\eqref{eq: negative Killing}, which in this $n=2$ case is given by $(X,Y) := - \kappa(X,Y) := - 4 \mathrm{Tr}(XY)$ for any $X,Y \in \mathfrak{su}(2)$. We also have 
\begin{equation}\label{eq:commutators}
[x,y] = \frac{1}{\sqrt{2}} t   \quad \textup{and} \quad [t,x] = \frac{1}{\sqrt{2}}y \quad \textup{ and } \quad [y,t] = \frac{1}{\sqrt{2}} x. 
\end{equation}

To compute $E_s$ as given in~\eqref{Es definition} we begin by observing that, since $\mathcal{B}$ is orthonormal, the map $\phi$ can be explicitly described as $\phi(x^*) = x, \quad \phi(y^*) = y, \quad \textup{ and } \quad \phi(t^*) = t.$
Thus $\phi(\xi t^* + \eta_1 x^* + \eta_2 y^*) = \xi t + \eta_1 x + \eta_2 y$. Throughout what follows, we write elements $p \in AN$ as
\begin{equation}\label{eq: AN coordinates} 
p = \begin{bmatrix}
a & z\\
0 & a^{-1}
\end{bmatrix}
\end{equation}
where $a>0$. It will also be helpful to define  
\begin{equation}\label{eq: Delta def} 
\Delta := \sqrt{\xi^2 + \eta_1^2 + \eta_2^2}.
\end{equation}
We start with the following lemma.

\begin{lemma}\label{lemma: Es coordinates} 
	Let $\lambda = \xi t^* + \eta_1 x^* + \eta_2 y^* \in \mathfrak{su}(2)^*$ and $s >0$. The matrix entries of  $p = E_s(\lambda ) \in AN$ as with respect to the coordinates introduced in~\eqref{eq: AN coordinates} are
	\begin{equation}\label{eq: def a and z} 
	a = \left(\cosh(s\Delta/\sqrt{2}) + \xi \frac{\sinh(s\Delta/\sqrt{2})}{\Delta}\right)^{-1/2}, \quad \text{and} \quad z = 	a (i\eta_1-\eta_2)\frac{\sinh(s\Delta/\sqrt{2})}{\Delta}.
	\end{equation} 
\end{lemma}

\begin{proof}
	We first compute the image of $\lambda = \xi t^* + \eta_1 x^* + \eta_2 y^* \in \mathfrak{su}(2)^*$ under the composition, which we denote $j_s$, of the map $\phi$ with the exponential map $\mathfrak{su}(2)  \xrightarrow{\exp(2s\sqrt{-1} \cdot ) } P$. Since $\phi(\xi t^* + \eta_1 x^* + \eta^2 y^*) = \xi t + \eta_1 x + \eta_2 y$ we obtain 
	\begin{align*}
	j_s(\xi t^{*} + \eta_1 x^{*} + \eta_2 y^*) &= \exp \left(\frac{s}{\sqrt{2}} \begin{bmatrix}
	-\xi  & i\eta_1 -\eta_2\\
	-i\eta_1 -\eta_2 & \xi
	\end{bmatrix}\right)
	\\
	&= \begin{bmatrix}
	\cosh(s\Delta/\sqrt{2}) - \xi \frac{\sinh(s\Delta/\sqrt{2})}{\Delta} & (i\eta_1-\eta_2)\frac{\sinh(s\Delta/\sqrt{2})}{\Delta} \\
	(-i\eta_1 - \eta_2)\frac{\sinh(s\Delta/\sqrt{2})}{\Delta} & \cosh(s\Delta/\sqrt{2}) + \xi \frac{\sinh(s\Delta/\sqrt{2})}{\Delta}
	\end{bmatrix}. 
	\end{align*}
	Next, we compute the inverse of $f$. By definition~\eqref{eq: def f Kstar to P}, $f(p) = p \cdot \tau(p)$ where $\tau(p) = \bar{p}^T$. Thus, in terms of the coordinates \eqref{eq: AN coordinates}, 
	\begin{align*}
	f\left(\begin{bmatrix}
	a & z \\
	0 & a^{-1}
	\end{bmatrix}\right) &= \begin{bmatrix}
	a^2 + |z|^2 & a^{-1}z\\
	a^{-1}\overline{z} & a^{-2}
	\end{bmatrix}. 
	\end{align*}
We thus have the following matrix equality 
	\[
	\begin{bmatrix}
	a^2 + |z|^2 & a^{-1}z\\
	a^{-1}\overline{z} & a^{-2}
	\end{bmatrix} = \begin{bmatrix}
	\cosh(s\Delta/\sqrt{2}) - \xi \frac{\sinh(s\Delta/\sqrt{2})}{\Delta} & (i\eta_1-\eta_2)\frac{\sinh(s\Delta/\sqrt{2})}{\Delta} \\
	(-i\eta_1 - \eta_2)\frac{\sinh(s\Delta/\sqrt{2})}{\Delta} & \cosh(s\Delta/\sqrt{2}) + \xi \frac{\sinh(s\Delta/\sqrt{2})}{\Delta}
	\end{bmatrix}. 
	\]
	from which the claims of the lemma follow. 
\end{proof}

It will be useful below to have notation for the real and imaginary parts of the complex number $z$ appearing in Lemma~\ref{lemma: Es coordinates}. From~\eqref{eq: def a and z} it is immediate that if $z = u+\sqrt{-1} v$ then 
\begin{equation}\label{eq: def u and v} 
u = -a \eta_2 \sinh(s \Delta/\sqrt{2})/\Delta \quad \textup{ and } \quad v = a \eta_1 \sinh(s\Delta/\sqrt{2})/\Delta.
\end{equation}

Recall that $\mathfrak{su}(2)^*$ is equipped with the $Ad^*$-invariant inner product $(,)$ induced from the $Ad$-invariant inner product $(,) = -\kappa(,)$ on $\mathfrak{su}(2)$. Given a non-zero element $\lambda \in \mathfrak{su}(2)^*$, let $Pr_\lambda$ denote the orthogonal projection from $\mathfrak{su}(2)^*$ onto the line spanned $\lambda$. Let $Pr_\lambda^\perp$ denote orthogonal projection onto the plane orthogonal to $\lambda$.  Explicitly, for any element $\zeta \in \mathfrak{su}(2)^*$, 
\begin{equation}
	Pr_\lambda(\zeta) = \frac{(\zeta,\lambda)}{\Delta^2}\lambda, \quad \text{and} \quad Pr_\lambda^\perp(\zeta) = \zeta - \frac{(\zeta,\lambda)}{\Delta^2}\lambda.
\end{equation}
With respect to the canonical identification $\mathfrak{su}(2)^* = T_\lambda \mathfrak{su}(2)^*$, the image of $Pr_\lambda^\perp$ is identified with the tangent space $T_\lambda \mathcal{O}_\lambda$, where $\mathcal{O}_\lambda$.  In what follows it will be useful to write $Pr_\lambda^\perp(\zeta)= ad_X \lambda$ for some Lie algebra element $X \in \k$. The following lemma provides a formula for such an element $X$.

\begin{lemma}\label{lemma: projections} 
	Let $\lambda = \xi t^* + \eta_1 x^* + \eta_2 y^* \in \mathfrak{su}(2)^*$. Assume $\lambda \neq 0$. Then
	\begin{equation}
		Pr_\lambda^\perp(t^*) = ad_{X_t}^*(\lambda), \quad Pr_\lambda^\perp(x^*) = ad_{X_x}^*(\lambda), \quad Pr_\lambda^\perp(y^*) = ad_{X_y}^*(\lambda), 
	\end{equation}
	where
	\begin{equation}
		X_t = \frac{\sqrt{2}}{\Delta^2}(\eta_2 x - \eta_1 y), \quad X_x = \frac{\sqrt{2}}{\Delta^2}(\xi y - \eta_2 t ), \quad X_y = \frac{\sqrt{2}}{\Delta^2}(\eta_1 t - \xi x).
	\end{equation}
\end{lemma}

\begin{proof}
	We check the formula for $Pr_\lambda^\perp(t^*)$. The other cases are similar. We compute 
	\begin{equation}
		ad_{X_t}^*\lambda = \frac{\sqrt{2}}{\Delta^2}(\eta_2 ad_{x}^*\lambda - \eta_1 ad_{y}^*\lambda)  = \frac{(\eta_1^2 + \eta_2^2)t^*-\xi\eta_2y^* - \xi\eta_1 x^*}{\Delta^2} = t^* -  \frac{\xi}{\Delta^2}\lambda =  Pr_\lambda^\perp(t^*). \qedhere
	\end{equation}

\end{proof}

We next give an explicit formula for $L_{p^{-1}}(\pi_{AN})_p$ from Lemma~\ref{lemma: piAN take 1} in the special case under consideration (i.e. $\k = \mathfrak{su}(2)$).  This requires fixing dual orthonormal bases. We have already fixed the explicit orthonormal basis $\{t,x,y\}$ of $\k$.  This determines a basis of $\mathfrak{an}$ which is dual to $\{t,x,y\}$ with respect to the pairing $\langle-,-\rangle_s$ which we denote as $\{t^*,x^*,y^*\}$. Explicitly,
\begin{equation}\label{eq: dual basis in AN}
x^*  =  \frac{s}{\sqrt{2}} \begin{bmatrix} 0 & \sqrt{-1} \\ 0 & 0 \end{bmatrix}, \quad 
y^* = -\frac{s}{\sqrt{2}} \begin{bmatrix} 0 & 1 \\ 0 & 0 \end{bmatrix}, \quad \text{and} \quad  
t^* = \frac{s}{2 \sqrt{2}} \begin{bmatrix} -1 & 0 \\ 0 & 1 \end{bmatrix}. 
\end{equation} 
Although this notation is very convenient in what follows, it introduces a potential source of confusion which we now attempt to dispel. The pairing $\langle-,-\rangle_s$ fixes an identification of $\mathfrak{an}$ with the abstract dual vector space $\k^*$. \emph{With respect to this identification}, the basis for $\mathfrak{an}$ that we just introduced \emph{is} the dual basis $\{t^*,x^*,y^*\}$ in $\k^*$. Despite this, the reader may  find it helpful to distinguish between the two spaces and thus also between their respective bases (for instance, if they are not accustomed to thinking about the moment map for a Hamiltonian $(SU(2),0)$ action as taking values in a vector space of upper triangular matrices).  In this case, the reader should note that the context will always dictate whether $t^*$, $x^*$, and $y^*$ are best understood as elements of $\k^*$ or $\mathfrak{an}$. For instance, in the following lemma $t^*,x^*,y^*$ are best understood as elements of $\mathfrak{su}(2)^*$ in the first sentence, and they are best understood as elements of $\mathfrak{an}$ in \eqref{Poisson}. 

\begin{lemma}\label{lemma: piAN take 2} 
	Let $\lambda=\xi t^* + \eta_1 x^* + \eta_2 y^* \in \mathfrak{su}(2)^*$ and let $p  = E_s(\lambda) = (a,u,v) \in AN$. Then the Poisson structure $\pi_{AN}$ on $SU(2)^*\cong AN$ satisfies 
	\begin{equation}\label{Poisson}
	L_{p^{-1}}(\pi_{AN})_p = \frac{a^{-1}u}{s} x^* \wedge t^* + \frac{a^{-1}v}{s} y^* \wedge t^* + \frac{a^{-2}}{2s} \left( u^2 + v^2 - a^2 + a^{-2}\right) x^* \wedge y^*.
	\end{equation}

\end{lemma}

In the course of the proof of Lemma~\ref{lemma: piAN take 2} we need the projection $Pr_{\mathfrak{an}}: \mathfrak{sl}(2,\C) \to \mathfrak{an}$ associated to the direct sum decomposition $\mathfrak{sl}(2,\C) \cong \mathfrak{su}(2,\C) \oplus \mathfrak{an}$. A straightforward computation yields the formula 
\begin{equation}\label{eq: proj to an}
Pr_{\mathfrak{an}}\left( \begin{bmatrix} \alpha & \beta \\ \gamma & - \alpha \end{bmatrix} \right) = 
\begin{bmatrix} Re(\alpha) & \beta + \overline{\gamma} \\ 0 & - Re(\alpha) \end{bmatrix}. 
\end{equation} 

\begin{proof}[Proof of Lemma~\ref{lemma: piAN take 2}]
	We begin by computing the terms $Pr_{\mathfrak{an}}(Ad_{p^{-1}}(b_i))$ where $b_i$ is set equal to $t, x$ and $y$ respectively. Let $p \in AN$ so $p = \begin{bmatrix} a& z \\ 0 & a^{-1} \end{bmatrix}$ for some $a>0, a \in \R$, and $z = u+ \sqrt{-1}v$. 
From~\eqref{eq: proj to an} it follows straightforwardly that 
	\[
	Pr_{\mathfrak{an}}(Ad_{p^{-1}}(t)) = \frac{\sqrt{-1}}{ \sqrt{2}} \begin{bmatrix} 0 & a^{-1}z \\ 0 & 0 \end{bmatrix}
	= \frac{a^{-1}u}{s} x^* +  \frac{a^{-1}v}{s} y^*
	\]
	and similar computations for $Pr_{\mathfrak{an}}(Ad_{p^{-1}}(x))$ and $Pr_{\mathfrak{an}}(Ad_{p^{-1}}(y))$ yield 
	\[
	Pr_{\mathfrak{an}}(Ad_{p^{-1}}(x)) 
		= -\frac{au}{s} t^*  + \frac{uv}{s} x^* - \frac{(a^{-2} - a^2 + u^2 - v^2)}{2s} y^* 
	\]
	and 
	\begin{equation*}
	Pr_{\mathfrak{an}}(Ad_{p^{-1}}(y))  = 
	 -\frac{av}{s} t^*  + \frac{(a^{-2} - a^2 - u^2 + v^2)}{2s} x^* - \frac{uv}{s} y^*. 
	\end{equation*}
	It is also straightforward to compute directly that 
	\begin{equation*} 
	Ad_{p^{-1}}(t^*)= t^* - a^{-1} vx^{*} + a^{-1} uy^{*} \, \, \textup{ and } \, \, 
Ad_{p^{-1}}(x^*) = a^{-2}x^* \, \, \textup{ and } \, \, 
	Ad_{p^{-1}}(y^*) = a^{-2}y^*. 
	\end{equation*} 
	Now by using the expression~\eqref{PoissonBivectorAN} of Lemma~\ref{lemma: piAN take 1}, the claim follows. 
\end{proof}

With the above lemmas in hand, we can now proceed to compute $\mathcal{D}(\omega_{can}) := \omega_{can} + \mu_L^* \Omega^s$ as a $2$-form on $T^*K \cong K \times \k^*$. We first address the special case of points of the form $(e,\lambda)$.  We choose an ordered frame $\{t, x, y\} \cup \{t^*, x^*, y^*\}$ for $T(T^*K)$ over the points $(e,\lambda)$, where $\{t,x,y\} =\mathcal{B}$ is the (ordered) basis of $\k = T_eK$ and $\{t^*, x^*, y^*\} = \mathcal{B}^*$ is the (ordered) basis of $\k^* \cong T_\lambda \k^*$. With respect to this frame, which restricts at each point $(e,\lambda)$ to a basis of $T(T^*K)$, the $2$-form $\mathcal{D}(\omega_{can})$ at $(e,\lambda)$ can be specified by a $6 \times 6$ skew-symmetric matrix, expressed in the block form 
\begin{equation}\label{eq: delin matrix} 
\begin{pmatrix} 
A & B \\ 
- B^t & C 
\end{pmatrix} 
\end{equation} 
where $A,B,C$ are $3 \times 3$ (real) matrices and $A$ and $C$ are skew-symmetric. Here we think of $A$ as the matrix recording the pairings of $\mathcal{B}$ with $\mathcal{B}$, while $B$ records the pairings of $\mathcal{B}$ against $\mathcal{B}^*$, and $C$ is the matrix of pairings $\mathcal{B}^*$ against $\mathcal{B}^*$. 
In what follows, we compute each of the matrices $A, B$ and $C$, in turn.

Lemma~\ref{lemma: case 1} below computes the matrix $A$.  Note that by skew-symmetry, in order to specify $A$, it suffices to compute the $3$ pairings given in the lemma. 

\begin{lemma}\label{lemma: case 1} Let $(e,\lambda) \in K\times \k^*$ and let $p  = E_s(\lambda) = (a,u,v) \in AN$.  Consider the basis elements $\mathcal{B} = \{x,y,t\}$ of $\mathfrak{su}(2)$ as elements of $T_{(e,\lambda)}(K \times \k^*) \cong T_eK \times T_\lambda k^* \cong \k \times \k^*$. 
Under this identification, we have 
	\[
		\mathcal{D}(\omega_{can})_{(e,\lambda)}(x, y) = \frac{a^{-2}}{2s} \left( u^2 + v^2 - a^2 + a^{-2}\right), \quad \mathcal{D}(\omega_{can})_{(e,\lambda)}(y,t) = \frac{a^{-1}v}{s}, \quad \mathcal{D}(\omega_{can})_{(e,\lambda)}(x,t) = \frac{a^{-1}u}{s}.
	\]
\end{lemma}

\begin{proof}
Let $p = E_s(\lambda)$ and $X,Y \in \{x,y,t\}$. 
Then $X = X_{T^*K}(e,\lambda)$ and $Y= Y_{T^*K}(e,\lambda)$. By Lemma~\ref{background useful fact}
\[
	\mathcal{D}(\omega_{can})_{(e,\lambda)}(X,Y) = (L_{p^{-1}})_*(\pi_{AN,s})_p(X,Y).
\]
On the right side, $X$ denotes  $\langle -,X\rangle_s$ in $\mathfrak{an}^*$ and similarly for $Y$. The result now follows from Lemma \ref{lemma: piAN take 2}. 
\end{proof}

\begin{remark} 
Notice that as $s \to 0$, the expressions in Lemma~\ref{lemma: case 1} limit as follows:  
\[
\lim_{s \to 0} \mathcal{D}(\omega_{can})_{(e,\lambda)}(x, y)  = \frac{\xi}{\sqrt{2}} , \quad 
\lim_{s \to 0} \mathcal{D}(\omega_{can})_{(e,\lambda)}(y,t) = \frac{\eta_1}{\sqrt{2}}, \quad 
\lim_{s \to 0}  \mathcal{D}(\omega_{can})_{(e,\lambda)}(t,x) =  \frac{\eta_2}{\sqrt{2}} 
\]
which matches $\omega_{can}$ as expected, as can be seen from a comparison with Example~\ref{omega can explicit}. 
\end{remark}

To compute the matrices $B$ and $C$, we need to compute the pushforward $d(E_s)_{\lambda}(\lambda)$ of the tangent vector $\lambda$ (considered as an element of $T_\lambda k^* \cong k^*$) at the point $\lambda \in \k^*$. This leads us to our next lemma. 

\begin{lemma}\label{lemma: dEs} 
Let $\lambda = \xi t^* + \eta_1 x^* + \eta_2 y^* \in \k^*$. 
For a real parameter $s$ define $\varepsilon := \sinh(s\Delta/\sqrt{2}) \Delta + \cosh(s\Delta/\sqrt{2}) \xi$. 
Then 
\begin{equation}\label{eq: Eslambda formula}
\Theta^L d(E_s)_\lambda(\lambda) = 
a^2 (\varepsilon t^* + \eta_1 x^* +  \eta_2 y^*). 
\end{equation} 
\end{lemma}

\begin{proof} 
Consider the path $\gamma(t) = (1+t)\lambda$ which goes through $\lambda$ at $t=0$
with $\frac{d}{dt} \big\vert_{t=0} \gamma(t) = \lambda \in T_\lambda k^*$. Thus $d(E_s)_\lambda(\lambda) = \frac{d}{dt} \big\vert_{t = 0} \begin{bmatrix} a(\gamma(t))  & z(\gamma(t)) \\ 0  & a^{-1}(\gamma(t)) \end{bmatrix}$. For simplicity we let $a(t) := a(\gamma(t))$ and $z(t) := z(\gamma(t))$, and $a(0) = a$. Note that $\Theta^L d(E_s)_\lambda(\lambda) = E_s(\lambda)^{-1} d(E_s)_\lambda(\lambda)$. Then 
\[
d(E_s)_\lambda(\lambda) = \begin{bmatrix} a'(0) & z'(0) \\ 0 & - a^{-2} \cdot a'(0) \end{bmatrix}
\quad \textup{ and } \quad 
E_s(\lambda)^{-1} d(E_s)_\lambda(\lambda) = \begin{bmatrix} a^{-1} a'(0) & a^{-1} z'(0) + z a^{-2} a'(0) \\ 0 & -a^{-1} a'(0) \end{bmatrix}. 
\]
Now straightforward calculations from~\eqref{eq: def a and z} yield 
$a'(0) = - \frac{s}{2\sqrt{2}} a^3 \varepsilon$
and
\[
z'(0) = a'(0) \left( \frac{(\sqrt{-1} \eta_1 - \eta_2)}{\Delta} \sinh\left( \frac{s\Delta}{\sqrt{2}}\right) \right) + a \frac{(\sqrt{-1}\eta_1 - \eta_2)}{\Delta} \cosh\left( \frac{s\Delta}{\sqrt{2}}\right) \frac{s \Delta}{\sqrt{2}}
\]
from which we
can compute 
\begin{equation}\label{eq: upper right corner} 
\begin{split} 
a^{-1} z'(0) + z a^{-2} a'(0) & = \frac{s}{\sqrt{2}} (\eta_2 - \sqrt{-1} \eta_1)\left( \frac{a^2 \varepsilon}{\Delta} \sinh(s \Delta/\sqrt{2})-\cosh(s\Delta/\sqrt{2})\right). 
\end{split} 
\end{equation} 
The claim now follows from the above and the definitions of $x^*, y^*, t^*$ in~\eqref{eq: dual basis in AN}, after simplifications using the definitions of $a$ and $\varepsilon$ and the identity $\cosh^2(\alpha) - \sinh^2(\alpha) = 1$ for a real number $\alpha$. 
\end{proof}

We now compute the matrix $B$ in~\eqref{eq: delin matrix}. 
By the moment map equation we have
\begin{equation}\label{eq: case 2 prelim} 
\begin{split}
(\mathcal{D}(\omega_{can}))_{(e,\lambda)}((X_{T^*K})_{(e,\lambda)}, (0, \eta)) & = \langle \Theta^L(d(E_s)_{\lambda}(\eta)), X \rangle_{E_s(\lambda)} . \\ 
\end{split} 
\end{equation} 
In order to evaluate the expression on the right, we write $\eta = Pr_\lambda(\eta) + Pr_\lambda^\perp(\eta)$ and evaluate $\langle \Theta^L(d(E_s)_{\lambda}(Pr_\lambda(\eta))), X \rangle_{E_s(\lambda)} $  using Lemma~\ref{lemma: dEs}. We evaluate $\langle \Theta^L(d(E_s)_{\lambda}(Pr_\lambda^\perp(\eta))), X \rangle_{E_s(\lambda)} $ using Lemma~\ref{lemma: projections}, equivariance of the map $E_s$, and the identity
\begin{equation}
\begin{split}
\langle \Theta^L(\mathcal{D}_Y(E_s)_{\lambda}(\lambda)), X \rangle_{E_s(\lambda)}  = (L_{E_s(\lambda)^{-1}})_*(\pi_{AN,s})_{E_s(\lambda)}(X,Y) \\ 
\end{split} 
\end{equation} 
which follows by the moment map equation for the dressing action of $K$ on the dressing orbit through $E_s(\lambda)^{-1}$.

The $3 \times 3$ matrix $B$ is given in the table below, where the elements in $\mathcal{B}^*$ label the columns and 
the elements in $\mathcal{B}$ label the rows.  
For the sake of this table we use the notation $\delta := u^2 + v^2 - a^2 + a^{-2}$.

\bigskip

\bgroup
\def\arraystretch{2.3}
\begin{tabular}{c|c|c|c}  
  & $t^*$ & $x^*$ & $y^*$ \\   \hline 
 
$ t$ &           $\frac{1}{\Delta^2} \left( \xi a^2 \varepsilon + \frac{\sqrt{2}a^{-1}}{ s}    (\eta_1 v - \eta_2 u) \right) $   &     $  \frac{1}{\Delta^2} \left( \eta_1 a^2 \varepsilon  - \frac{\sqrt{2} a^{-1} \xi v}{s} \right) $       &         $ \frac{1}{\Delta^2} \left( \eta_2 a^2 \varepsilon + \frac{\sqrt{2} a^{-1} \xi u}{s} \right)    $   \\ \hline 

 $x$ &     $ \frac{1}{\Delta^2} \left( \xi a^2 \eta_1 - \frac{ \eta_1 a^{-2} \delta}{\sqrt{2} s} \right)  $  &      
 $    \frac{1}{\Delta^2} \left( a^2 \eta_1^2 + \frac{\xi a^{-2} \delta}{\sqrt{2} s} - \frac{\sqrt{2} a^{-1} \eta_2 u}{s} \right) $ 
      &        $    \frac{1}{\Delta^2} \left( \eta_2 a^2 \eta_1 + \frac{\sqrt{2} a^{-1} \eta_1 u}{s} \right) $ 
          \\   \hline 
 
 $y$ &     $  \frac{1}{\Delta^2} \left( \xi a^2 \eta_2 - \frac{\eta_2 a^{-2} \delta}{\sqrt{2} s} \right) $          &     
 $ \frac{1}{\Delta^2} \left( \eta_1 a^2 \eta_2 - \frac{\sqrt{2}  a^{-1} \eta_2 v}{s} \right) $ 
      &       $ \frac{1}{\Delta^2} \left( a^2 \eta_2^2 + \frac{\sqrt{2} a^{-1} \eta_1 v}{s} + \frac{\xi a^{-2} \delta}{\sqrt{2} s} \right) $ 
\end{tabular} 
\egroup

\bigskip

\begin{remark} 
Taking the limit of the entries as $s$ approaches $0$, we can see that $\mathcal{D}(\omega_{can})$ approaches $\omega_{can}$ as expected. Consider, for instance, the entry in the above table corresponding to $(t, t^*)$. Since $\lim_{s \to 0} a^2 \varepsilon = \xi$, $\lim_{s \to 0} a^{-1} v/s = \eta_1/\sqrt{2}$, and $\lim_{s \to 0} a^{-1}u/s = - \eta_2/\sqrt{2}$, we conclude that
\[
\lim_{s \to 0} \frac{\xi}{\Delta^2} a^2 \varepsilon + \frac{\sqrt{2}a^{-1}}{\Delta^2 s}    (\eta_1 v - \eta_2 u) = \frac{\xi^2}{\Delta^2} + 
\frac{\eta_1^2}{\Delta^2} + \frac{\eta_2^2}{\Delta^2} = \frac{\Delta^2}{\Delta^2} = 1.
\]
This agrees with the coefficient of $(t^*)^L\wedge t$ in the expression for $\omega_{can}$ given in Example~\ref{omega can explicit}.  The limits of the other entries in the table above can be computed similarly. One obtains in the limit as $s\to 0$ the identity $3 \times 3$ matrix as desired. 
\end{remark}

Lastly, we compute the matrix $C$ in~\eqref{eq: delin matrix}.  As for the matrix $A$, by the skew-symmetry of $C$ it suffices to compute $3$ of the entries of $C$. Note first that $(\omega_{can})_{(e,\lambda)}((0,\eta_1), (0, \eta_2)) = 0$ for any $\eta_1, \eta_2$ due to properties of the canonical symplectic structure on any cotangent bundle. Therefore we have $$(\mathcal{D}(\omega_{can}))_{(e,\lambda)}((0, \eta_1), (0, \eta_2)) = (\Omega^s)_\lambda(\eta_1, \eta_2)$$ where we have applied Lemma~\ref{lemma: dmu}. To compute the RHS of this equation, we decompose $\eta_i = Pr_\lambda(\eta_i) + Pr^{\perp}_\lambda(\eta_i)$  for $i=1,2$ and use that 
$(\Omega^s)_\lambda(Pr_\lambda(\eta_1), Pr_\lambda(\eta_2)) = 0$ since $\Omega^s$ is skew-symmetric and $Pr_\lambda(\eta_1)$ and $Pr_\lambda(\eta_2)$ are parallel. Equation~\eqref{eq: Omega s def} then allows us to compute explicitly the pairings for $\eta_i \in \mathcal{B}^*$. 
In the course of the computations it is useful to notice that 
$\delta = 2 a^2 \varepsilon \frac{\sinh(s\Delta/\sqrt{2})}{\Delta}$ and therefore that $\eta_1 \delta = 2 a \varepsilon v$ and $\eta_2 \delta = -2a \varepsilon u$. The details are left to the reader and we record the results in the lemma below.

\begin{lemma}\label{lemma: case 3} 
Following the notation above, we have 

\begin{equation*}
(\mathcal{D}(\omega_{can}))_{(e,\lambda)}((0, x^*), (0, y^*)) = \frac{\sqrt{2}\xi}{\Delta^4}\left( \frac{\sinh(s\Delta/\sqrt{2})}{s\Delta/\sqrt{2}} \left( \xi(\xi-\varepsilon) + \Delta^2 \right) - \left( a^2\xi(\xi-\varepsilon) + \Delta^2 \right) \right) + \frac{\sqrt{2}a^2}{\Delta^2}(\xi-\varepsilon), 
\end{equation*} 

\begin{equation*} 
(\mathcal{D}(\omega_{can}))_{(e,\lambda)}((0, t^*), (0, x^*)) =   \frac{\sqrt{2}\eta_2}{\Delta^4}\left( \frac{\sinh(s\Delta/\sqrt{2})}{s\Delta/\sqrt{2}} \left( \xi(\xi-\varepsilon) + \Delta^2 \right) - \left( a^2\xi(\xi-\varepsilon) + \Delta^2 \right) \right), 
\end{equation*} 

\begin{equation*} 
(\mathcal{D}(\omega_{can}))_{(e,\lambda)}((0,y^*),(0,t^*)) =   \frac{\sqrt{2} \eta_1}{\Delta^4} \left( \frac{\sinh(s\Delta/\sqrt{2})}{s\Delta/\sqrt{2}} \left( \xi(\xi-\varepsilon) + \Delta^2 \right) - \left( a^2\xi(\xi-\varepsilon) + \Delta^2 \right) \right). 
\end{equation*} 

\end{lemma}

\bigskip

\bigskip

In the computations recorded above, we restricted above to the special case of points of the form $(e,\lambda)$. Finally, to complete the discussion, we briefly sketch how $\mathcal{D}(\omega_{can})$ may be computed at arbitrary points. The basis of the strategy is the simple observation that by Lemma~\ref{lemma: delinearized Hamiltonian by H}, the $2$-form $\mathcal{D}(\omega_{can})$ is invariant under the \emph{right} multiplication action of $K$ on $T^*K$, i.e. for any $h \in K$ we have $R_h^*(\mathcal{D}(\omega_{can})) = \mathcal{D}(\omega_{can})$.  Specifically, an arbitrary  point $(k,\nu)$ can be brought to a point $(e, Ad^*_k(\nu)=\lambda)$ by acting by $R_k$. We have already expressed $\mathcal{D}(\omega_{can})$ at $(e, Ad^*_k(\nu))$ as 
\[
(\mathcal{D}(\omega_{can}))_{(e,Ad^*_k(\nu))} = \frac{1}{2}\sum_{a,b\in \mathcal{B}} \alpha_{a,b} a^* \wedge b^* + 
\sum_{a,b\in \mathcal{B}} \beta_{a,b} a^* \wedge b  +  \frac{1}{2}
\sum_{a, b \in \mathcal{B}} \gamma_{a,b} a \wedge b 
\]
where the coefficients $\alpha_{a,b}, \beta_{a,b}, \gamma_{a,b}$ are as computed above. Then in terms of the basis of right-invariant $1$-forms on $K$ (for the first factor of $T^*K \cong K \times \k^*$) and the basis $\mathcal{B} = (\mathcal{B}^*)^*$ for the second factor we have 
\[
(\mathcal{D}(\omega_{can}))_{(k, \nu)} = \frac{1}{2}\sum_{a,b\in \mathcal{B}} \alpha_{a,b} R^*_h a^* \wedge R^*_h b^* + 
\sum_{a,b\in \mathcal{B}} \beta_{a,b} R^*_h a^* \wedge Ad_{k^{-1}}(b)  + \frac{1}{2}
\sum_{a, b \in \mathcal{B}} \gamma_{a,b} Ad_{k^{-1}}(a) \wedge Ad_{k^{-1}}(b)
\]
and a computation of $Ad_{k^{-1}}(\eta)$ for $\eta \in \mathcal{B}$ then yields an explicit expression for $\mathcal{D}(\omega_{can})$ in terms of the invariant $1$-forms and $\mathcal{B}$.




\begin{thebibliography}{9}

\bibitem{Alekseev} A.~Alekseev, \emph{On Poisson actions of compact Lie groups on symplectic manifolds}, J.~Differential Geom.~\textbf{47} (1997), pp.~241--256.

\bibitem{AlekseevHoffmanLaneLi} A.~Alekseev, B.~Hoffman, J.~Lane, and Y.~Li, \emph{Action-angle coordinates on coadjoint orbits and multiplicity free spaces from partial tropicalization}, arxiv:2003.13621.

\bibitem{AlekseevMeinrekenWoodward} A.~Alekseev, E.~Meinrenken, and C.~Woodward, \emph{Linearization of Poisson actions and singular values of matrix products}, Ann.~Inst.~Fourier (Grenoble) \textbf{51} (2001) no.~6, pp.~1691--1717.

\bibitem{CannasDaSilva} 
A.~C.~da Silva, \emph{Lectures on Symplectic Geometry}, Springer-Verlag, Berlin, Heidelberg, 2008.

\bibitem{GuilleminSternberg} V.~Guillemin and S.~Sternberg, \emph{A normal form for the moment map}, In S.~Sternberg, editor, \emph{Differential geometric methods in mathematical physics (Jerusalem, 1982)}, vol.~6, pp.~161--175. Reidel, 1984.

\bibitem{PoissonStructures} 
C.~Laurent-Gengoux, A.~Pichereau, and P.~Vanhaecke, \emph{Poisson Structures}, Springer-Verlag, Berlin, Heidelberg, 2013.

\bibitem{LuMomentMap} J.~H.~Lu, \emph{Momentum mappings and reduction of {P}oisson actions}, In P.~Dazord and A.~Weinstein, editors, \emph{Symplectic Geometry, Groupoids, and Integrable Systems}, pp.209--226. Springer, 1991.

\bibitem{Lu} J.~H.~Lu, \emph{Classical dynamical $r$-matrices and homogeneous Poisson structures on $G/H$ and $K/T$}, Commun.~Math.~Phys. \textbf{212} (2000), pp.~337--370.

\bibitem{LuThesis} J.~H.~Lu, \emph{Multiplicative and affine Poisson structures on Lie groups}, PhD Thesis, University of California at Berkeley. 1990.

\bibitem{LuWeinstein} J.~H.~Lu and A.~Weinstein. \emph{Poisson Lie groups, dressing transformations, and Bruhat decompositions}, J.~Differential Geom.~\textbf{31} (1990), pp.~501--526.

\bibitem{Marle} C.-M.~Marle, \emph{Mod\`ele d'action Hamiltonienne d'un groupe de Lie sur une vari\'et\'e symplectique}, Rend.~Sem.~Mat.~Univ.~Politec.~Torino, \textbf{43} (2014), pp.~187--205.







\end{thebibliography}
\end{document}